\pdfoutput=1

\newif\ifTwoColumn
\newif\ifTechReport

\TechReporttrue  

\ifTwoColumn
        \documentclass[10pt, journal, twoside, doublecolumn]{IEEEtran}
        \usepackage{etex}

        \IEEEoverridecommandlockouts

\else
        \ifTechReport
                \documentclass[11pt,onecolumn,a4paper]{IEEEtran}
                \usepackage{etex}
                \usepackage[left = 2.5cm, right = 2.5cm, top = 2.5cm, bottom = 2.5cm]{geometry}
                \usepackage{amsthm}
                \usepackage{parskip}
                \makeatletter
                \def\thm@space@setup{%
                  \thm@preskip=\parskip \thm@postskip=0pt
                }
                \makeatother

        \else
                \documentclass[12pt,draftcls,onecolumn,a4paper]{IEEEtran}
                \usepackage{amsthm}
                \IEEEoverridecommandlockouts
        \fi
\fi

\usepackage{cite}
\usepackage{hyperref}

\ifCLASSINFOpdf
   \usepackage[pdftex]{graphicx}
\else
\fi
\usepackage[cmex10]{amsmath}
\usepackage{amssymb}
\usepackage{algorithm}
\usepackage{algpseudocode}

\usepackage{array}

\usepackage[caption=false,font=footnotesize]{subfig}
\usepackage{mathtools}  
\usepackage{booktabs}   
\usepackage{units}      
\usepackage{multirow}   
\usepackage{xcolor,calc}        
\usepackage{balance}
\usepackage{tikz}       
\usetikzlibrary{shapes,arrows, matrix, calc}
\usetikzlibrary{external} 
\usepackage{pgfplots}
\usepackage{setspace}
\usepackage{adjustbox} 

\definecolor{redCol}{rgb}{0.6350, 0.0780, 0.1840}
\definecolor{greenCol}{rgb}{0.4660,  0.6740,    0.1880}
\definecolor{blueCol}{rgb}{0,   0.4470,   0.7410}

\DeclareMathOperator{\Tr}{Tr}

\begin{document}
%


\title{High-Speed Finite Control Set Model Predictive Control for Power Electronics}

%
%
%

\author{Bartolomeo~Stellato,~\IEEEmembership{Student Member,~IEEE,}  Tobias~Geyer,~\IEEEmembership{Senior Member,~IEEE} and~Paul~J.~Goulart,~\IEEEmembership{Member,~IEEE}
\thanks{This work was supported by the People Programme (Marie Curie Actions) of the European Union’s Seventh Framework Programme (FP7/2007-2013) under REA grant agreement no 607957 (TEMPO).}
\thanks{B.\ Stellato and P.J.\ Goulart are with the University of Oxford, Parks Road, Oxford, OX1 3PJ, U.K. (e-mail: \mbox{bartolomeo.stellato@eng.ox.ac.uk}; \mbox{paul.goulart@eng.ox.ac.uk}).}
\thanks{T.\ Geyer is with the ABB Corporate Research Centre, 5405 Baden-Dattwil, Switzerland (e-mail: \mbox{t.geyer@ieee.org}).}
}
\maketitle

\begin{abstract}
Common approaches for direct model predictive control (MPC) for current reference tracking in power electronics suffer from the high computational complexity encountered when solving integer optimal control problems over long prediction horizons. We propose an efficient alternative method based on approximate dynamic programming, greatly reducing the computational burden and enabling sampling times below $\unit[25]{\mu s}$. Our approach is based on the offline estimation of an infinite horizon value function which is then utilized as the tail cost of an MPC problem. This allows us to reduce the controller horizon to a very small number of stages while simultaneously improving the overall controller performance. Our proposed algorithm was implemented on a small size FPGA and validated on a variable speed drive system with a three-level voltage source converter. Time measurements showed that our algorithm requires only $\unit[5.76]{\mu s}$ for horizon $N=1$ and  $\unit[17.27]{\mu s}$ for $N=2$, in both cases outperforming state of the art approaches with much longer horizons in terms of currents distortion and switching frequency. To the authors' knowledge, this is the first time direct MPC for current control has been implemented on an FPGA solving the integer optimization problem in real-time and achieving comparable performance to formulations with long prediction horizons.
\end{abstract}

\begin{IEEEkeywords}
Approximate dynamic programming, value function approximation, drive systems,  finite control set, model predictive control (MPC).
\end{IEEEkeywords}

%
\IEEEpeerreviewmaketitle

\ifTwoColumn \else
  \ifTechReport \else
    \newpage
  \fi
\fi


\section{Introduction}
\label{sec:introduction}

\IEEEPARstart{A}{mong} the control strategies adopted in power electronics, model predictive control (MPC)~\cite{Rawlings:2014wt} has recently gained popularity due to its various advantages~\cite{Cortes:2008cv}. MPC has been shown to outperform traditional control methods mainly because of its ease in handling time-domain constraint specifications that can be imposed by formulating the control problem as a constrained optimization problem. Due to its structure, MPC can be applied to a variety of power electronics topologies and operating conditions providing a higher degree of flexibility than traditional approaches.

With the recent advances in convex optimization techniques~\cite{nocedal2006numerical}, it has been possible to apply MPC to very fast constrained linear systems with continuous inputs by solving convex quadratic optimization problems within microseconds~\cite{Jerez:2014jw}. However, when dealing with nonlinear systems or with integer inputs, the optimal control problems are no longer convex and it is harder to find optimal solutions. Sequential quadratic programming (SQP)~\cite{nocedal2006numerical} has gained popularity because of its ease in iteratively approximating nonconvex continuous control problems as convex quadratic programs, which can be solved efficiently. Moreover, integrated perturbation analysis (IPA) has been recently combined with SQP methods (IPA-SQP)~\cite{Ghaemi:2009di} and applied to power electronics~\cite{Xie:2012ca} by solving the approximated quadratic optimization problem at a given time instant using a perturbed version of the problem at the previous time instant, thereby reducing the number of iterations required at each time step. However, there are still two orders of magnitude difference in achievable computation time compared to results obtained for linear systems~\cite{Jerez:2014jw} and further advances are required to apply these methods to very fast dynamical systems.

In power electronics, many conventional control strategies applied in industry are based on proportional-plus-integral (PI) controllers providing continuous input signals to a modulator that manages conversion to discrete switch positions. Direct MPC~\cite{Geyer:2005us} instead combines the current control and the modulation problem into a single computational problem, providing a powerful alternative to conventional PI controllers. With direct MPC, the manipulated variables are the switch positions, which lie in a discrete and finite set, giving rise to a switched system. Therefore, this approach does not require a modulator and is often referred to as \emph{finite control set} MPC.

Since the manipulated variables are restricted to be integer-valued, the optimization problem underlying direct MPC is provably $\mathcal{NP}$-hard~\cite{bertsimas2005optimization}.
In power electronics these optimization problems are often solved by complete enumeration of all the possible solutions, which grow exponentially with the prediction horizon~\cite{Quevedo:2004db}.
Since long horizons are required to ensure stability and good closed-loop performance~\cite{Morari:1999ho},  direct MPC quickly becomes intractable for real-time applications. As a consequence, in cases when reference tracking of the converter currents is considered, the controller horizon is often restricted to one~\cite{Cortes:2008cv}.
Recently, the {IPA-SQP} method has been applied to a finite control set MPC~\cite{Nademi:2016wn} to efficiently approximate the optimization problem in the case of nonlinear systems. However, no particular attention is paid to reducing the number of integer combinations that must be evaluated, which is  at the source of the most significant computational issues.
Despite attempts to overcome the computational burden of integer programs in power electronics~\cite{Karamanakos:2014kh}, the problem of implementing these algorithms on embedded systems remains open.

A recent technique introduced in~\cite{Geyer:ij} and benchmarked in~\cite{Geyer:fva} reduces the computational burden of direct MPC when increasing the prediction horizon. In that work the optimization problem was formulated as an integer least-squares (ILS) problem and solved using a tailored branch-and-bound algorithm, described as \emph{sphere decoding}~\cite{Hassibi:2005gu}, generating the optimal switching sequence at each time step. Although this approach appears promising relative to previous work, the computation time required to perform the sphere decoding algorithm for long horizons (i.e.\ $N=10$), is still far slower than the sampling time typically required, i.e.~${T_s = \unit[25]{\mu s}}$.
In the very recent literature, some approaches have been studied to improve the computational efficiency of the sphere decoding algorithm. In particular, in~\cite{PetrosKaramanakos:2014we} a method based on a lattice reduction algorithm decreased the average computational burden of the sphere decoding. However, the worst case complexity of this new reformulation is still exponential in the horizon length. In~\cite{Karamanakos:2015jl}, heuristic search strategies for the sphere decoding algorithm are studied at the expense of suboptimal control performance. Even though a floating point complexity analysis of the algorithms is presented in these works, no execution times and no details about fixed-point implementation are provided. Furthermore, there currently exists no embedded implementation of a direct MPC algorithm for current control achieving comparable performance to formulations with long prediction horizons.

This paper introduces a different method to deal with the direct MPC problem. In contrast to common formulations~\cite{Kouro:2009gs} where the switching frequency is controlled indirectly via penalization of the input switches over the controller horizon, in this work the system dynamics are augmented to directly estimate the switching frequency. 
Our approach allows the designer to set the desired switching frequency a priori by penalizing its deviations from this estimate. Thus, the cost function tuning can be performed more easily than with the approach in~\cite{Geyer:fva} and~\cite{Geyer:ij}, where a tuning parameter spans the whole frequency range with no intuitive connection to the desired frequency. To address the computational issues of long prediction horizons, we formulate the tracking problem as a regulation one by augmenting the state dynamics and cast it in the framework of approximate dynamic programming~(ADP)~\cite{bertsekas1996dynamic}. The infinite horizon value function is approximated using the approach in~\cite{deFarias:2003wq} and~\cite{Wang:2014em} by solving a semidefinite program (SDP)~\cite{Vandenberghe:1996fy} offline. This enables us to shorten the controller horizon by applying the estimated tail cost to the last stage to maintain good control performance. In~\cite{2016arXiv160207273B} the authors applied a similar approach to stochastic systems with continuous inputs, denoting the control law as the ``iterated greedy policy''.

As a case study, our proposed approach is applied to a variable-speed drive system consisting of a three-level neutral point clamped voltage source inverter connected to a medium-voltage induction machine. The plant is modelled as a linear system with a switched three-phase input with equal switching steps for all phases.

Closed loop simulations in MATLAB in steady state operation showed that with our method even very short prediction horizons give better performance than the approach in~\cite{Geyer:fva} and~\cite{Geyer:ij} with much longer planning horizons.

We have implemented our algorithm on a small size Xilinx Zynq FPGA (xc7z020) in fixed-point arithmetic and verified its performance with hardware in the loop (HIL) tests of both steady-state and transient performance. The results achieve almost identical performance to closed-loop simulations and very fast computation times, allowing us to comfortably run our controller  within the $\unit[25]{\mu s}$ sampling time.

The remainder of the paper is organized as follows. In Section~\ref{sec:DRIVE} we describe the drive system case study and derive the physical model. In Section~\ref{sec:dmpc} the direct MPC problem is derived by augmenting the state dynamics and approximating the infinite horizon tail cost using ADP. 
 Section~\ref{sec:Framework for Performance Evaluation} describes all of the physical parameters of the model used to verify our approach. In Section~\ref{sec:Achievable Performance in Steady State} we present closed-loop simulation results on the derived model in steady state operation to characterize the achievable performance of our method. In Section~\ref{sec:FPGA Implementation} we describe the hardware setup, the algorithm and all the FPGA implementation details. In Section~\ref{sec:Hardware In The Loop Tests} HIL tests are performed in steady-state and transient operation. Finally, we provide conclusions in Section~\ref{sec:conclusions}.

In this work we use normalized quantities by adopting the per unit system (pu). The time scale $t$ is also normalized using the base angular velocity $\omega_b$ that in this case is $\unitfrac[2 \pi 50]{rad}{s}$, i.e.\ one time unit in the per unit system corresponds to $\unit[1/\omega_b]{s}$. 
Variables in the three-phase system (\textit{abc}) are denoted $\boldsymbol{\xi}_{abc} = \begin{bmatrix} \xi_a & \xi_b & \xi_c\end{bmatrix}^\top$. It is common practice to transform phase variables to $\boldsymbol{\xi}_{\alpha\beta}$ in the stationary orthogonal $\alpha\beta$ coordinates by $\boldsymbol{\xi}_{\alpha\beta} = \boldsymbol{P}\boldsymbol{\xi}_{abc}$.
The inverse operation can be  performed as $\boldsymbol{\xi}_{abc}= \boldsymbol{P}^{\dagger}\boldsymbol{\xi}_{\alpha\beta}$. The matrices $\boldsymbol{P}$ and $\boldsymbol{P}^{\dagger}$ are the Clarke transform and its pseudoinverse respectively, i.e.\
\begin{equation*}
	\boldsymbol{P} = \dfrac{\strut 2}{\strut 3}\begin{bmatrix}
 1 & -\dfrac{\strut 1}{\strut 2} & - \dfrac{\strut 1}{\strut 2}\\
 0 & \dfrac{\strut \sqrt{3}}{\strut 2} & -\dfrac{\strut \sqrt{3}}{\strut 2}
 \end{bmatrix},
 \quad \boldsymbol{P}^{\dagger} = \begin{bmatrix} 1 & 0\\
	-\dfrac{\strut 1}{\strut 2} & \dfrac{\strut \sqrt{3}}{\strut 2}\\
	-\dfrac{\strut 1}{\strut 2} & -\dfrac{\strut \sqrt{3}}{\strut 2}
 \end{bmatrix}.
\end{equation*}



\section{Drive System Case Study}\label{sec:DRIVE}
\begin{figure}[hbt]
\begin{center}
  \ifTwoColumn
  \includegraphics[width = \columnwidth]{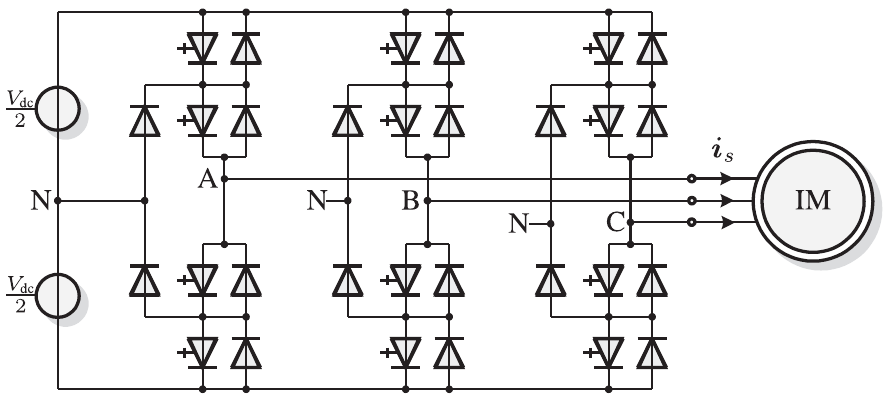}
  \else
  \includegraphics[width = 0.6\columnwidth]{img/inv_ind_machine}
  \fi
  \caption{Three-level three-phase neutral point clamped (NPC) voltage source inverter driving an induction motor with a fixed neutral point potential.}
  \label{img:drive:complete_phys_system}
  \end{center}
\end{figure}

In this work we consider a variable-speed drive system as shown in Figure~\ref{img:drive:complete_phys_system}, consisting of a three-level neutral point clamped (NPC) voltage source inverter driving a medium-voltage (MV) induction machine. The total dc-link voltage $V_{dc}$ is assumed constant and the neutral point potential $N$ fixed.

In most modern approaches to control variable-speed drive systems, the control is split into two cascaded loops. The outer loop controls the machine speed by manipulating the torque reference. The inner loop controls the torque and the fluxes by manipulating the voltages applied to the stator windings of the machine. Our approach focuses on the inner loop. The reference torque is converted into stator currents references that must be tracked and the controller manipulates the stator voltages by applying different inverter switch positions.

\subsection{Physical Model of the Inverter}
The switch positions in the three phase legs can be described by the integer input variables $u_a, u_b, u_c \in\left\{-1, 0, 1\right\}$, leading to phase voltages $\{-\frac{V_{dc}}{2}, 0, \frac{V_{dc}}{2}\}$, respectively. Hence, the output voltage of the inverter is given by
\begin{equation}\label{eq:drive:output_voltages}
	\boldsymbol{v}_{\alpha\beta} = \frac{V_{dc}}{2}\boldsymbol{u}_{\alpha\beta} = \frac{V_{dc}}{2}\boldsymbol{P}\boldsymbol{u}_{sw},
\end{equation}
where $\boldsymbol{u}_{sw} = \begin{bmatrix} u_a & u_b & u_c\end{bmatrix}^\top$.


\subsection{Physical Model of the Machine}
Hereafter we derive the state-space model of the squirrel-cage induction machine in the $\alpha\beta$ plane. Since we are considering a current control problem, it is convenient to use the stator current $\boldsymbol{i}_{s,\alpha\beta}$ and the rotor flux $\boldsymbol{\psi}_{r,\alpha\beta}$ as state variables. The model input is the stator voltage $\boldsymbol{v}_{s,\alpha\beta}$ which is equal to the inverter output voltage in~\eqref{eq:drive:output_voltages}. The model parameters are: the stator and rotor resistances $R_s$ and $R_r$; the mutual, stator and rotor reactances $X_{m}, X_{ls}$ and $X_{lr}$, respectively; the inertia $J$; and the mechanical load torque $T_l$. Given these quantities, the continuous-time state equations\cite{Krause:2002vy,Holtz:G982etC8} are

\begin{equation}\label{eq:machine:physical_model}
\begin{aligned}
\frac{\mathrm{d}\boldsymbol{i}_{s}}{\mathrm{d}t} & = - \frac{1}{\tau_s} \boldsymbol{i}_{s} + \left(\frac{1}{\tau_r}\boldsymbol{I} - \omega_r \begin{bmatrix}
 0 & -1\\ 1 & 0
 \end{bmatrix}
\right)\frac{X_m}{D}\boldsymbol{\psi}_r + \frac{X_r}{D}\boldsymbol{v}_s\\
\frac{\mathrm{d}\boldsymbol{\psi}_{r}}{\mathrm{d}t} & = \frac{X_m}{\tau_r} \boldsymbol{i}_{s} - \frac{1}{\tau_r}\boldsymbol{\psi}_{r} + \omega_r \begin{bmatrix}
 0 & -1\\ 1 & 0
 \end{bmatrix} \boldsymbol{\psi}_r\\
\frac{\mathrm{d}\omega_{r}}{\mathrm{d}t} & = \frac{1}{J}\left(T_e - T_l\right),
\end{aligned}
\end{equation}
where $D \coloneqq X_s X_r - X_m^2$ with $X_s \coloneqq X_{ls} + X_m$ and $X_r \coloneqq X_{lr} + X_m$, and $\boldsymbol{I}$ represents the $2 \times 2$ identity matrix. To simplify the notation, we have dropped the subscripts $\alpha\beta$ from all vectors in~\eqref{eq:machine:physical_model}. Moreover, $\tau_s \coloneqq X_r D/\left(R_sX_r^2 + R_r X_m^2\right)$ and $\tau_r \coloneqq X_r /R_r$ are the transient stator and the rotor time constants respectively. The electromagnetic torque is given by
\begin{equation}\label{eq:drive_system:torque}
  T \coloneqq \frac{X_m}{X_r} \left(\boldsymbol{\psi}_r \times \boldsymbol{i}_s\right).
\end{equation}
The rotor speed $\omega_r$ is assumed to be constant within the prediction horizon.
For prediction horizons in the order of a few milliseconds this is a mild assumption.
\subsection{Complete Model of the Physical System}\label{sec:physic_system:states}
Given the models of the drive and of the induction motor in \eqref{eq:drive:output_voltages} and \eqref{eq:machine:physical_model} respectively, the state-space model in the continuous time domain can be described as
\begin{subequations}\label{eq:complete_phys:continuous_time}
\begin{align}
\frac{\mathrm{d}\boldsymbol{x}_{ph}(t)}{\mathrm{d}t} & = \boldsymbol{D}\boldsymbol{x}_{ph}(t) + \boldsymbol{E}\boldsymbol{u}_{sw}(t)\label{eq:complete_phys:continuous_time_a}\\
\boldsymbol{y}_{ph}(t)& = \boldsymbol{F}\boldsymbol{x}_{ph}(t),
\end{align}
\end{subequations}
where the state vector $\boldsymbol{x}_{ph} = \begin{bmatrix}i_{s, \alpha} & i_{s, \beta} & \psi_{r,\alpha} & \psi_{r,\beta}\end{bmatrix}^\top$ includes the stator current and rotor flux in the $\alpha\beta$ reference frame. The output vector is taken as the stator current, i.e. $\boldsymbol{y}_{ph} = \boldsymbol{i}_{s, \alpha\beta}$. The matrices $\boldsymbol{D}, \boldsymbol{E}$ and $\boldsymbol{F}$ are defined in Appendix~\ref{app:phys_sys_mat}.

The state-space model of the drive can be converted into the discrete-time domain using exact Euler discretization. By integrating~\eqref{eq:complete_phys:continuous_time_a} from $t = k\hat{T}_s$ to $t = (k+1)\hat{T}_s$ and keeping $\boldsymbol{u}_{sw}(t)$ constant during each interval and equal to $\boldsymbol{u}_{sw}(k)$, the discrete-time model becomes
\begin{subequations}\label{eq:complete_phys:discrete_time}
\begin{align}
\boldsymbol{x}_{ph}(k+1) & = \boldsymbol{A}_{ph}\boldsymbol{x}_{ph}(k) + \boldsymbol{B}_{ph}\boldsymbol{u}_{sw}(k)\\
\boldsymbol{y}_{ph}(k)& = \boldsymbol{C}_{ph}\boldsymbol{x}_{ph}(k),
\end{align}
\end{subequations}
with matrices $\boldsymbol{A}_{ph} \coloneqq e^{\boldsymbol{D}\hat{T}_s}$, ${\boldsymbol{B}_{ph} \coloneqq - \boldsymbol{D}^{-1}\left(\boldsymbol{I} - \boldsymbol{A}_{ph}\right)\boldsymbol{E}}$, $\boldsymbol{C}_{ph} \coloneqq \boldsymbol{F}$ and $k\in \mathbb{N}$. $\boldsymbol{I}$ is an identity matrix of appropriate dimensions. Although the sampling time is $T_s = \unit[25]{\mu s}$, we use the discretization interval $\hat{T}_s = T_s \omega_b$ for consistency with our per unit system, where $\omega_b$ is the base frequency.



\section{Model Predictive Current Control}\label{sec:dmpc}

\subsection{Problem Description}
\label{sub:Problem Description}

Our control scheme must address two conflicting objectives simultaneously.
On the one hand, the distortion of the stator currents $\boldsymbol{i}_{s}$ cause iron and copper losses in the machine leading to thermal losses. Because of the limited cooling capability of the electrical machine, the stator current distortions have to be kept as low as possible.
On the other hand, high frequency switching of the inputs $\boldsymbol{u}_{sw}$ produces high power losses and stress on the semiconductor devices.
Owing to the limited cooling capability in the inverter, we therefore should minimize the switching frequency of the integer inputs.


Note that the effect of the inverter switchings on the torque ripples can be improved during the machine design. In particular, increasing the time constants of the stator and the rotor $\tau_s$ and $\tau_r$ can reduce the amplitude of the torque ripples by decreasing the derivative of the currents $\boldsymbol{i}_s$ and fluxes $\boldsymbol{\varphi}_r$. This is achieved naturally when dealing with machines with higher power.
Thanks to the flexibility of model based controller designs such as MPC, different machine dynamics influencing the torque ripples are automatically taken into account by the controller, which adapts the optimal inputs computation depending on the plant parameters. Thus, any improvements during the machine design can be optimally exploited by adapting the internal model dynamics in the controller.
Another similar approach is to include LCL filters between the inverter and the motor to decrease the high frequency components of the currents; see~\cite{Scoltock:2015dl}. These approaches allow operation at lower switching frequencies with low THD at the same time. However, it is sometimes impossible to change the machine's physical configuration and it is necessary to operate the inverter at high switching frequencies to satisfy high performance requirements in terms of stator currents distortion. For all these reasons there is an unavoidable tradeoff between these two criteria.

The controller sampling time plays an important role in the distortion and switching frequency tradeoff. Depending on the precision required in defining the inverter switching times, the controller is discretized with higher (e.g. $\unit[125]{\mu s}$) or lower (e.g. $\unit[25]{\mu s}$) sampling times. Higher sampling times define a more coarse discretization grid leading to less precise definition of the switching instants, but more available time to perform the computations during the closed-loop cycles. Lower sampling times, on the other hand, lead to improved controller accuracy while reducing the allowed computing time. However, for the same switching frequency, longer sampling times produce higher distortions. Ideally, the sampling time should be chosen as low as possible to have the highest possible accuracy.

In contrast to the common approaches in direct MPC where the switching frequency is minimized, in this work we penalize its difference from the desired frequency which is denoted by~$f_{sw}^*$. This is motivated by the fact that inverters are usually designed to operate at a specific nominal switching frequency.

The current distortion is measured via the total harmonic distortion (THD). Given an infinitely long time-domain current signal $i$ and its fundamental component $i^*$ of constant magnitude, the THD is proportional to the root mean square (RMS) value of the difference $i - i^*$. Hence, we can write for one phase current
\begin{equation}\label{eq:dmpc:thd}
	THD \sim \lim_{M\rightarrow \infty} \sqrt{\frac{1}{M}\sum_{k=0}^{M-1}(i(k) - i^*(k))^2},
\end{equation}
with $M\in \mathbb{N}$. For the three-phase current $\boldsymbol{i}_{abc}$ and its reference $\boldsymbol{i}_{abc}^*$ the THD is proportional to the mean value of~\eqref{eq:dmpc:thd} over the phases. It is of course not possible to calculate the THD in real time within our controller computations because of finite storage constraints.

The switching frequency of the inverter can be identified by computing
the average frequency of each active semiconductor device. As
displayed in Figure~\ref{img:drive:complete_phys_system}, the total number of switches in all three phases is 12, and for each switching transition by one step up or down in a phase one semiconductor device is turned on.
Thus, the number of \emph{on} transitions occurring between time step $k-1$ and $k$ is given by the $1$-norm of the difference of the inputs vectors: $\|\boldsymbol{u}_{sw}(k) - \boldsymbol{u}_{sw}(k-1)\|_1$.

Given a time interval centered at the current time step $k$ from $k-M$ to $k+M$, it is possible to estimate the switching frequency by counting the number of \textit{on} transitions over the time interval and dividing the sum by the interval's length $2MT_s$. We then can average over all the semiconductor switches by dividing the computed fraction by~$12$.
At time $k$, the switching frequency estimate can be written as
\begin{equation}\label{eq:dmpc:fswM}
	f_{sw,M}(k) \coloneqq \frac{1}{12\cdot 2M T_s}\sum_{i = -M}^{M} \|\boldsymbol{u}_{sw}(k + i) - \boldsymbol{u}_{sw}(k + i - 1)\|_1,
\end{equation}
which corresponds to a non-causal finite impulse response (FIR) filter of order $2M$.
The true average switching frequency is the limit of this quantity as the window length goes to infinity
\begin{equation}\label{eq:dmpc:fsw}
	f_{sw} \coloneqq \lim_{M\to \infty}f_{sw, M}(k),
\end{equation}
and does not depend on time $k$.

The $f_{sw}$ computation brings similar issues as the THD. In addition to finite storage constraints, the part of the sum regarding the future signals produces a non-causal filter that is impossible to implement in a real-time control scheme.


These issues in computing THD and $f_{sw}$ are addressed in the following two sections via augmentation of our state space model to include suitable approximation schemes for both quantities.

\subsection{Total Harmonic Distortion}\label{sec:dmpc:thd}
According to~\eqref{eq:dmpc:thd}, the THD in the three-phase current is proportional to the mean value of $(i_{s,a} - i_{s,a}^*)^2 + (i_{s,b} - i_{s,b}^*)^2 + (i_{s,c} - i_{s,c}^*)^2$. As shown in~\cite{Geyer:ij}, the THD is also proportional to the stator current ripple in the $\alpha\beta$ coordinate system, i.e.
\begin{equation}
	THD \sim \lim_{M\rightarrow \infty} \sum_{k=0}^{M-1}\|\boldsymbol{e_i}(k)\|_2^2,
\end{equation}
where we have introduced the error signal ${\boldsymbol{e}_{\boldsymbol{i}}(k) \coloneqq  \boldsymbol{i}_{s, \alpha\beta}(k) - \boldsymbol{i}_{s, \alpha\beta}^{*}(k)}$. It is straightforward to show~\cite{5060554} that the stator current reference during steady-state operation at rated frequency is given by
\begin{equation}\label{eq:dmpc:curr_alphabeta}
	\boldsymbol{i}_{s, \alpha\beta}^{*}(k) = \begin{bmatrix}\sin\left(k\right) & -\cos\left(k\right)\end{bmatrix}^{\top}.
\end{equation}
Hence, in order to minimize the THD, we minimize the squared $2$-norm of vector $\boldsymbol{e_i}$ over all future time steps.
%
We also introduce a discount factor $\gamma\in (0,1)$ to normalize the summation preventing it from going to infinity due to persistent tracking errors. The cost function related to THD minimization is therefore
\begin{equation}\label{eq:dmpc:thd_cost_fun}
	\sum_{k=0}^{\infty} \gamma^k\left\| \boldsymbol{e}_{\boldsymbol{i}}(k) \right\|_2^2.
\end{equation}

In order to construct a regulation problem, we include the oscillating currents from~\eqref{eq:dmpc:curr_alphabeta} as two additional uncontrollable states $\boldsymbol{x}_{osc}=\boldsymbol{i}_{s, \alpha\beta}^{*}$ within our model of the system dynamics. The ripple signal $\boldsymbol{e}_{\boldsymbol{i}}(k)$ is then modeled as an output defined by the difference between two pairs of system states.

\subsection{Switching Frequency}\label{sec:dmpc:frequency_filter}
To overcome the difficulty of dealing with the filter in~\eqref{eq:dmpc:fswM}, we consider only the past input sequence, with negative time shift giving a causal FIR filter estimating $f_{sw}$. This filter is approximated with an infinite impulse response (IIR) one whose dynamics can be modeled as a linear time invariant (LTI) system. Note that future input sequences in~\eqref{eq:dmpc:fswM} are taken into account inside the controller prediction.

Let us define three binary phase inputs denoting whether each phase switching position changed at time $k$ or not, i.e.\
\begin{equation}\label{eq:dmpc:rel_p_u1}
	\boldsymbol{p}(k) \coloneqq \begin{bmatrix} p_a(k) & p_b(k) & p_c(k) \end{bmatrix}^\top\in \{0,1\}^3,
\end{equation}
with
\begin{equation}\label{eq:dmpc:rel_p_u2}
	p_s(k) = \left\| u_s(k) - u_s(k-1)\right\|_1, \; s \in \{a, b, c\}.
\end{equation}
It is straightforward to show that the following second order IIR filter will approximate the one-sided version of the FIR filter in~\eqref{eq:dmpc:fswM}~\cite{Oppenheim:1997tq}:
\begin{align}\label{eq:sw_filter_dynamics}
	\boldsymbol{x}_{flt}(k+1) &= \underbrace{\begin{bmatrix}
 a_1 & 0\\ 1 - a_1 & a_2
 \end{bmatrix}}_{\boldsymbol{A_{flt}}}
\boldsymbol{x}_{flt}(k) + \underbrace{\dfrac{\strut 1 - a_2}{\strut 12 T_s}\begin{bmatrix}1 & 1 & 1\\ 0 & 0 & 0\end{bmatrix}}_{\boldsymbol{B}_{flt}}\boldsymbol{p}(k)\\
	\hat{f}_{sw}(k) &= \begin{bmatrix}
 0 & 1
 \end{bmatrix}
\boldsymbol{x}_{flt}(k),
	\end{align}
where $\hat{f}_{sw}(k)$ is the estimated switching frequency and $\boldsymbol{x}_{flt}(k)$ is the filter state.
The two poles at $a_1 = 1 - 1/r_1$ and $a_2=1-1/r_2$ with $r_1, r_2 >> 0$ can be tuned to shape the behavior of the filter.  Increasing $a_1,a_2$ make the estimate smoother, while decreasing $a_1,a_2$ gives a faster estimation with more noisy values.

We denote the difference between the approximation $\hat{f}_{sw}(k)$ and the target frequency $f_{sw}^*$ by $e_{sw}(k) := \hat{f}_{sw}(k) - f^{*}(k)$. Therefore, the quantity to be minimized in order to bring the switching frequency estimate as close to the target as possible is
\begin{equation}\label{eq:dmpc:sw_freq_cost_fun}
	 \sum_{k=0}^{\infty}\delta\gamma^k\left\|e_{sw}(k)\right\|_2^2,
\end{equation}
where $\delta \in \mathbb{R}_{+}$ is a design  parameter included to reflect the importance of this part of the cost relative to the THD component.

Finally, we can augment the state space to include the filter dynamics
and the target frequency by adding the states $\begin{bmatrix}\boldsymbol{x}_{flt}^\top &
f_{sw}^{*}\end{bmatrix}^\top$ so that the control inputs try to drive
the difference between two states to zero. Since the physical states are expressed in the per unit (pu) system with values around $1$, in order to have these augmented states within the same order of magnitude we will normalize them by the desired frequency $f_{sw}^*$ defining $\boldsymbol{x}_{sw} = \begin{bmatrix}
(1/f_{sw}^{*})\boldsymbol{x}_{flt}^\top & 1
\end{bmatrix}$ and the matrices $\boldsymbol{A}_{sw} = \text{blkdiag}(\boldsymbol{A}_{sw}, 1)$, $\boldsymbol{B}_{sw} = \begin{bmatrix}\boldsymbol{B}_{flt}^\top & \boldsymbol{0}_{1\times 3}^\top \end{bmatrix}^\top$.

\subsection{MPC Problem Formulation}
\label{sub:MPC Problem Formulation}

Let us define the complete augmented state as 
\begin{equation}\label{eq:dmpc:state_def}
	\boldsymbol{x}(k) \coloneqq \begin{bmatrix}\boldsymbol{x}_{ph}(k)^\top & \boldsymbol{x}_{osc}(k)^\top &  \boldsymbol{x}_{sw}(k)^\top & \boldsymbol{u}_{sw}(k-1)^\top \end{bmatrix}^\top,
\end{equation}
with ${}\boldsymbol{x}(k) \in \mathbb{R}^{9} \times \{-1, 0, 1\}^3$ and total state dimension $n_x = 12$. The vector $\boldsymbol{x}_{ph}$ represents the physical system from Section~\ref{sec:physic_system:states}, $\boldsymbol{x}_{osc}$ defines the oscillating states of the sinusoids to track introduced in Section~\ref{sec:dmpc:thd}, ${\boldsymbol{u}_{sw}(k-1)}$ are additional states used to keep track of the physical switch positions at the previous time step, and $\boldsymbol{x}_{sw}$ are the states related to the switching filter from Section~\ref{sec:dmpc:frequency_filter}.

The system inputs are defined as
\begin{equation*}
	\boldsymbol{u}(k) \coloneqq \begin{bmatrix}\boldsymbol{u}_{sw}(k)^\top & \boldsymbol{p}(k)^\top \end{bmatrix}^\top\in \mathbb{R}^{n_u},
\end{equation*}
where $\boldsymbol{u}_{sw}$ are the physical switch positions and $\boldsymbol{p}$ are the three binary inputs entering in the frequency filter from Section~\ref{sec:dmpc:frequency_filter}.  The input dimension is $n_u = 6$. To simplify the notation let us define the matrices $\boldsymbol{G}$ and $\boldsymbol{T}$ to obtain $\boldsymbol{u}_{sw}(k)$ and $\boldsymbol{p}(k)$ from $\boldsymbol{u}(k)$ respectively, i.e.\ such that ${\boldsymbol{u}_{sw}(k) = \boldsymbol{G}\boldsymbol{u}(k)}$ and ${\boldsymbol{p}(k) = \boldsymbol{T}\boldsymbol{u}(k)}$. Similarly, to obtain $\boldsymbol{u}_{sw}(k-1)$ from $\boldsymbol{x}(k)$ we define a matrix $\boldsymbol{W}$ so that $\boldsymbol{u}_{sw}(k-1) = \boldsymbol{W}\boldsymbol{x}(k)$.

The MPC problem with horizon $N\in \mathbb{N}$ can be written as%
\begin{subequations}\label{eq:compmpc:problem_full}
\begin{alignat}{2}
&\underset{\boldsymbol{u}(k)}{\text{minimize}} \quad && \sum_{k=0}^{N-1}\gamma^{k}\ell(\boldsymbol{x}(k))+ \gamma^{N}V(\boldsymbol{x}(N))\label{eq:compmpc:problem_full:cost_fun}\\
&\text{subject to} \quad  &&  \boldsymbol{x}(k+1)  = \boldsymbol{A}\boldsymbol{x}(k) + \boldsymbol{B}\boldsymbol{u}(k) \label{eq:compmpc:problem_full:dyn}\\
         &&& \boldsymbol{x}(0) = \boldsymbol{x}_0 \label{eq:compmpc:problem_full:initst}\\
          &&& \boldsymbol{x}(k) \in \mathcal{X},\;\boldsymbol{u}(k) \in \mathcal{U}(\boldsymbol{x}(k)),
 \end{alignat}
\end{subequations}
where the stage cost is defined combining the THD and the switching frequency penalties in~\eqref{eq:dmpc:thd_cost_fun} and~\eqref{eq:dmpc:sw_freq_cost_fun} respectively as
\begin{equation*}
	\ell(\boldsymbol{x}(k)) = \left\|\boldsymbol{C}\boldsymbol{x}(k)\right\|_2^2 = \left\|\boldsymbol{e}_{\boldsymbol{i}}(k)\right\|_2^2 + \delta \left\|e_{sw}(k)\right\|_2^2.
\end{equation*}
The tail cost $V(\boldsymbol{x}(N))$ is an approximation of the infinite horizon tail that we will compute in the next section using approximate dynamic programming (ADP).
The matrices $\boldsymbol{A}$, $\boldsymbol{B}$ and $\boldsymbol{C}$ define the extended system dynamics and the output vector; they can be derived directly from the physical model~\eqref{eq:complete_phys:discrete_time} and from the considerations in Sections~\ref{sec:dmpc:frequency_filter} and~\ref{sec:dmpc:thd}.

The input constraint set is defined as
\begin{subequations}\label{eq:dmpc:constraints_u}
	\begin{align}
		\mathcal{U}(\boldsymbol{x}(k)) := &\{
		 \left\|\boldsymbol{T}\boldsymbol{u}(k)\right\|_\infty\leq 1\label{eq:dmpc:constraints_u:fsw2},\\
		&- \boldsymbol{T}\boldsymbol{u}(k) \leq \boldsymbol{G}\boldsymbol{u}(k) - \boldsymbol{W}\boldsymbol{x}(k)\leq \boldsymbol{T}\boldsymbol{u}(k),\label{eq:dmpc:constraints_u:fsw1}\\
& \boldsymbol{G}\boldsymbol{u}(k) \in \left\{-1, 0, 1\right\}^{3}\label{eq:dmpc:constraints_u:integr}\},
	\end{align}
\end{subequations}
where constraint~\eqref{eq:dmpc:constraints_u:fsw1} defines the relationship between $\boldsymbol{u}_{sw}$ and $\boldsymbol{p}$ from~\eqref{eq:dmpc:rel_p_u1} and~\eqref{eq:dmpc:rel_p_u2}. Constraint~\eqref{eq:dmpc:constraints_u:fsw2} together with~\eqref{eq:dmpc:constraints_u:fsw1} defines the switching constraints ${\left\| \boldsymbol{u}_{sw}(k) - \boldsymbol{u}_{sw}(k-1)\right\|_{\infty}\leq 1}$ imposed to avoid a shoot-through in the inverter positions that could damage the components. Finally,~\eqref{eq:dmpc:constraints_u:integr} enforces integrality of the switching positions.

It is straightforward to confirm that the number of switching sequence combinations grows exponentially with the horizon length $N$, i.e. $3^{3N} = 27^N$. The problem therefore becomes extremely difficult to solve for even modest horizon lengths.

Observe that the controller tuning parameters are $\delta$, which defines the relative importance of the THD and $f_{sw}$ components of the cost function, and $r_1,r_2$, which shape the switching frequency estimator.

\begin{figure}
\begin{center}
  \tikzstyle{block} = [draw, thick, fill=blueCol!20, rectangle, inner sep=5pt]
  \tikzstyle{sum} = [draw, circle, node distance=2.84em, inner sep=0pt]
  \tikzstyle{branch} = [circle,inner sep=0pt,minimum size=1mm,fill=black,draw=black]
  \tikzstyle{input} = [coordinate]
  \tikzstyle{output} = [coordinate]

  \begin{tikzpicture}[auto, node distance=5.68em, >=stealth', every node/.style={align=center}]

    \node [input] (torque) {};
    \node [block, right of = torque, node distance = 4.83em] (osc) {OSC};
    \node [block, right of = osc, node distance = 7.38em, minimum height=8.52em,yshift = -2.84em] (mpc) {MPC};
    \node [block, below of=mpc, node distance = 6.39em] (flt) {FLT};
    \node [block, below of=flt, node distance = 2.84em] (delay) {$z^{-1}$};
    \node [branch,right of=mpc, node distance=5.112em, yshift = 1.42em] (branchU) {};
    \node[coordinate, right of=mpc, node distance=4.26em, yshift = -1.42em](p){};
    \node[coordinate, left of=mpc, node distance=5.68em, yshift = -2.84em](pprev){};
    \node[coordinate, left of=mpc, node distance=7.1em, yshift = -0.9372em](uprev){};
    \node [block,right of=branchU, node distance=4.544em,minimum height=4.26em, minimum width = 4.26em] (motor) {MOTOR};
    \node [block, below of=motor, node distance = 8.52em] (observer) {OBS};
    \node[output, below of=observer, node distance = 1.42em] (xphys1){};
    \node[output, below of=delay, left of =xphys1] (xphys2){};
    \node[output, left of=mpc, node distance =11.36em, yshift = 0.9372em] (xphys3){};

    \draw [->] (torque) -- node [pos = 0.4] {$T^*(k)$}(osc);
    \draw [->] (osc) -- node [pos=0.55] {$\boldsymbol{x}_{osc}(k)$} (osc-|mpc.west);
    \draw [->] (branchU-|mpc.east) -- node [pos=0.55] {$\boldsymbol{u}_{sw}(k)$} (branchU) -- (motor);
    \draw [->] (motor.280) -- node {$\boldsymbol{i}_s(k)$}(observer.north-|motor.280);
    \draw [->] (motor.260) -- node [anchor = east]{$\omega_r(k)$}(observer.north-|motor.260);
    \draw [->] (mpc.east|-p) -- node [pos=0.5] {$\boldsymbol{p}(k)$} (p) |- (flt);
    \draw [->] (flt) -|  (pprev) -- node [pos=0.54] {$\boldsymbol{x}_{sw}(k)$}(pprev-|mpc.west);
    \draw [->] (branchU) |- (delay);
    \draw [->] (delay) -| (uprev) -- node[pos = 0.5] {$\boldsymbol{u}_{sw}(k-1)$}(uprev-|mpc.west);
    \draw [->] (observer.south) -- (xphys1) |- node [pos =0.1]{$\boldsymbol{x}_{ph}(k)$}(xphys2) -| (xphys3) -- node [pos = 0.81]{$\boldsymbol{x}_{ph}(k)$}(xphys3-|mpc.west);
		\draw[thick, dotted] ($(mpc.north) + (0.0,0.5)$) -| ($(branchU) + (0.25,0.0)$) |- ($(delay.south) + (0.0, -0.5)$) -| node[pos = 0.12, anchor = south east]{CONTROLLER}($(osc.west) + (-0.25,0.0)$) |- ($(mpc.north) + (0.0,0.5)$);
  \end{tikzpicture}
\end{center}
\caption{
Block diagram of the control loop. The controller within the dotted line receives the desired torque $T^*(k)$ and the motor states $\boldsymbol{x}_{ph}(k)$ providing the switch position $\boldsymbol{u}_{sw}(k)$. 
}
\label{fig:block_diagram}
\end{figure}
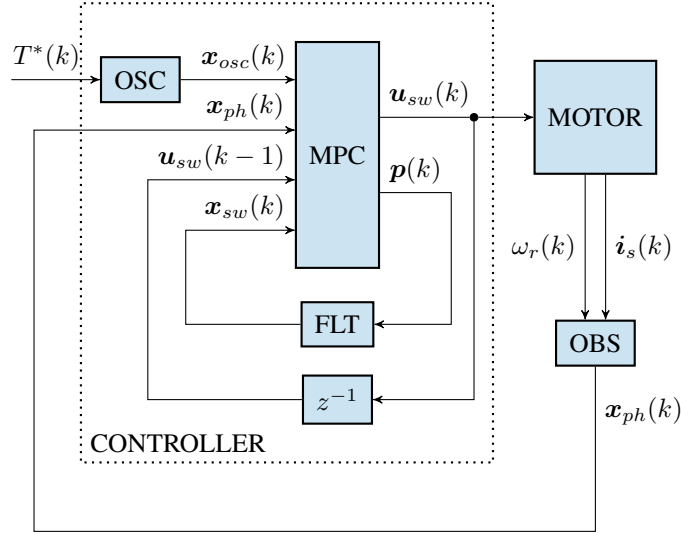

\subsection{Control Loop}
\label{sub:Control Loop}

The complete block diagram is shown in Figure~\ref{fig:block_diagram}. The desired torque $T^*$ determines the currents $\boldsymbol{x}_{osc}$ by setting the initial states of the oscillator OSC. The motor speed $\omega_r$ and the stator currents $\boldsymbol{i}_s$ are measured directly from the machine and used by the observer OBS providing the physical states of the motor $\boldsymbol{x}_{ph}$. The auxiliary inputs $\boldsymbol{p}$ are fed into the filter FLT estimating the switching frequency in $\boldsymbol{x}_{sw}$. The switch positions $\boldsymbol{u}_{sw}$ go through a one step delay and are exploited again by the MPC formulation.

Following a receding horizon control strategy, at each stage~$k$ the problem~\eqref{eq:compmpc:problem_full} is solved, obtaining the optimal sequence $\left\{\boldsymbol{u}^{\star}(k)\right\}_{k=0}^{N-1}$ from which only $\boldsymbol{u}^{\star}(0)$ is applied to the switches. At the next stage $k+1$, given new vectors $\boldsymbol{x}_{osc}(k), \boldsymbol{x}_{ph}(k), \boldsymbol{u}_{sw}(k-1)$ and $\boldsymbol{x}_{sw}(k)$ as in Figure~\ref{fig:block_diagram} a new optimization problem is then solved providing an updated optimal switching sequence, and so on. The whole control algorithm, appearing within the dotted line, runs within $\unit[25]{\mu s}$.

\subsection{Approximate Dynamic Programming}

The goal of this section is to compute a value function approximation $V^{adp}$ for an infinite horizon version of~\eqref{eq:compmpc:problem_full}. The function $V^{adp}$ is used as a tail cost in~\eqref{eq:compmpc:problem_full}.

Let $V^*(\boldsymbol{z})$ be the value function evaluated in $\boldsymbol{z}$, i.e.\ the optimal value of the objective of our control problem starting at state $\boldsymbol{z}$
\begin{equation*}
	V^*(\boldsymbol{z}) = \min_{\boldsymbol{u}\in \mathcal{U}(\boldsymbol{z})}\left\{\sum_{k=0}^{\infty}\gamma^{k}l(\boldsymbol{x}(k), \boldsymbol{u}(k))\right\},
\end{equation*}
subject to the system dynamics~\eqref{eq:compmpc:problem_full:dyn}.
For notational convenience, we will drop the time index $k$ from the vectors in this section.
The main idea behind dynamic programming is that the function $V^*$ is the unique solution to the equation
\begin{equation*}
	V^*(\boldsymbol{z}) = \min_{\boldsymbol{u}\in \mathcal{U}(\boldsymbol{z})}\left\{l(\boldsymbol{z}, \boldsymbol{u}) + \gamma V^*\left(\boldsymbol{A}\boldsymbol{z} + \boldsymbol{B}\boldsymbol{u}\right)\right\} \quad \forall \boldsymbol{z},
\end{equation*}
known as the Bellman equation. The right-hand side can be written as monotonic operator $\mathcal{T}$ on $V^*$, usually referred to as the Bellman operator: $V^* = \mathcal{T}V^*$.
Once $V^*$ is known, the optimal control policy for our problem starting at state $\boldsymbol{z}$ can be found as
\begin{equation*}
	\boldsymbol{\psi}^*(\boldsymbol{z}) = \arg\min_{\boldsymbol{u}\in \mathcal{U}(\boldsymbol{z})}\left\{l(\boldsymbol{z}, \boldsymbol{u}) + \gamma V^*\left(\boldsymbol{A}\boldsymbol{z} + \boldsymbol{B}\boldsymbol{u}\right)\right\},
\end{equation*}
subject to constraints~\eqref{eq:compmpc:problem_full:dyn}.

Unfortunately, solutions to the Bellman equation can only be solved analytically in a limited number of special cases; e.g.\ when the state and inputs have small dimensions or when the system is linear, unconstrained and the cost function is quadratic~\cite{Kalman:1964ha}. For more complicated problems, dynamic programming is limited by the so-called curse of dimensionality; storage and computation requirements tend to grow exponentially with the problem dimensions. Because of the integer switches in the power converter analyzed in this work, it is intractable to compute the optimal infinite horizon cost and policy and, hence, systematic methods for approximating the optimal value function  offline are needed.

Approximate dynamic programming~\cite{bertsekas1996dynamic} consists of various techniques for estimating $V^*$ using knowledge from the system dynamics, fitted data through machine learning or iterative learning through simulations.

\subsubsection*{Approximation via Iterated Bellman Inequalities} the approach developed in~\cite{deFarias:2003wq} and~\cite{Wang:2014em} relaxes the Bellman equation into an inequality
\begin{equation}\label{eq:adp:bellman_inequality}
	V^{adp}(\boldsymbol{z}) \leq \min_{\boldsymbol{u}\in \mathcal{U}(\boldsymbol{z})}\left\{l(\boldsymbol{z}, \boldsymbol{u}) + \gamma V^{adp} \left(\boldsymbol{A}\boldsymbol{z} + \boldsymbol{B}\boldsymbol{u}\right)\right\}, \quad \forall \boldsymbol{z},
\end{equation}
or, equivalently, using the Bellman operator: $V^{adp} \leq
\mathcal{T}V^{adp}$.

The set of functions $V^{adp}$ that satisfy the Bellman inequality are underestimators of the optimal value function $V^*$. This happens because, if $V^{adp}$ satisfies $V^{adp}\leq \mathcal{T}V^{adp}$, then by the monotonicity of the operator $\mathcal{T}$ and value iteration convergence, we can write
\begin{equation*}
	V^{adp} \leq \mathcal{T} V^{adp} \leq \mathcal{T}\left(\mathcal{T}V^{adp}\right) \leq \dots \leq \lim_{i\to\infty}\mathcal{T}^{i}V^{adp} = V^{*}.
\end{equation*}
The Bellman inequality is therefore a sufficient condition for underestimation of $V^*$. In~\cite{Wang:2014em} the authors show that by iterating inequality~\eqref{eq:adp:bellman_inequality}, the conservatism of the approximation can be reduced. The iterated Bellman inequality is defined as:
\begin{equation*}
	V^{adp}(\boldsymbol{z}) \leq \mathcal{T}^{M}V^{adp},
\end{equation*}
where $M>1$ is an integer defining the number of iterations. This inequality is equivalent to the existence of functions $V^{adp}_{i}$ such that
\begin{equation*}
	V^{adp} \leq \mathcal{T}V^{adp}_{1}, \quad V^{adp}_{1}\leq \mathcal{T}V^{adp}_2, \quad \dots V^{adp}_{M-1}\leq \mathcal{T}V^{adp}.
\end{equation*}
By defining $V^{adp}_0=V^{adp}_M=V^{adp}$, we can rewrite the iterated inequality as
\begin{equation}\label{eq:adp:iter_bellman_inequality2}
	V^{adp}_{i-1} \leq \mathcal{T}V^{adp}_{i}, \quad i = 1,\dots,M,
\end{equation}
where $V^{adp}_{i}$ are the iterates of the value function.

To make the problem tractable, we will restrict the iterates to the finite-dimensional subspace spanned by the basis functions $V^{(j)}$ defined in ~\cite{deFarias:2003wq},\cite{Wang:2014em}:
\begin{equation}\label{eq:adp:iter_bellman_inequality_subspace}
	V^{adp}_{i} = \sum_{j = 1}^{K}\alpha_{ij}V^{(j)}, \quad i = 0, \dots, M-1.
\end{equation}
The coefficients $\alpha_i$ will be computed by solving a Semidefinite Program (SDP)~\cite{Vandenberghe:1996fy}.

The rewritten iterated Bellman inequality in~\eqref{eq:adp:iter_bellman_inequality2} suggests the following optimization problem for finding the best underestimator for the value function $V^*$
\begin{subequations}\label{eq:adp:underestim_first}:
 \begin{alignat}{2}
&\text{maximize} \quad & & \int_{\mathcal{X}}V^{adp}(\boldsymbol{z})c(\mathrm{d}\boldsymbol{z})\\
&\text{subject to} \quad & & V^{adp}_{i-1}(\boldsymbol{z}) \!\leq\! \min_{\boldsymbol{u}\in \mathcal{U}(\boldsymbol{z})}\left\{l(\boldsymbol{z}, \boldsymbol{u}) + \gamma V^{adp}_{i}(\boldsymbol{A}\boldsymbol{z} + \boldsymbol{B}\boldsymbol{u})\right\}\label{eq:adp:iter_ineq_constr}\\
&&& \forall \boldsymbol{z} \in \mathbb{R}^{6}\times \{-1, 0, 1\},\quad i=1,\dots,M,\\
&&& V^{adp}_0=V^{adp}_M=V^{adp},
 \end{alignat}
\end{subequations}
where $c(\cdot)$ is a non-negative measure over the state space. On the chosen subspace~\eqref{eq:adp:iter_bellman_inequality_subspace}, the inequality~\eqref{eq:adp:iter_ineq_constr} is convex in the coefficients $\alpha_{ij}$. To see this, note that the left-hand side is affine in $\alpha_{ij}$. Moreover, for a fixed $\boldsymbol{u}$ the argument of $\min$ on the right-hand side is affine in $\alpha_{ij}$ while the $\min$ of affine functions is concave.

The solution to~\eqref{eq:adp:underestim_first} is the function spanned by the chosen basis that maximizes the $c$-weighted $1$-norm defined in the cost function while satisfying the iterated Bellman inequality~\cite{deFarias:2003wq}.
Hence, $c(\cdot)$ can be regarded as a distribution 
giving more importance to regions of the state space where we would like a better approximation.

Following the approach in~\cite{Wang:2014em}, we make use of quadratic candidate functions of the form
\begin{equation}\label{eq:adp:quad_fun}
	V^{adp}_i(\boldsymbol{z}) = \boldsymbol{z}^{\top}\boldsymbol{P}_i \boldsymbol{z} + 2 \boldsymbol{q}_i^\top \boldsymbol{z} + r_i, \quad i=0,\dots,M,
\end{equation}
where ${\boldsymbol{P}_i \in \mathbb{S}^{n_x},}\;{\boldsymbol{q}_i\in \mathbb{R}^{n_x},}\; {r_i \in \mathbb{R},}\; i = 0,\dots,M$.

If we denote $\boldsymbol{\mu}_c \in \mathbb{R}^{n_x}$ and $\boldsymbol{\Sigma}_c \in \mathbb{S}_{+}^{n_x}$ as the mean and the covariance matrix of measure $c(\cdot)$ respectively, by using candidate functions as in~\eqref{eq:adp:quad_fun} the cost function of problem~\eqref{eq:adp:underestim_first} becomes
\begin{equation*}
	\int_{\mathcal{X}}V^{adp}(\boldsymbol{z})c(\mathrm{d}\boldsymbol{z}) = \Tr\left(\boldsymbol{P}_0 \boldsymbol{\Sigma}_c\right) + 2\boldsymbol{q}_{0}^\top \boldsymbol{\mu}_c + r_0.
\end{equation*}
We now focus on rewriting the constraint~\eqref{eq:adp:iter_ineq_constr} as a Linear Matrix Inequality (LMI)~\cite{Boyd:1994tc}. We first remove the $\min$ on the right-hand side by imposing the constraint for every admissible $\boldsymbol{u} \in \mathcal{U}(\boldsymbol{x}_0)$ and obtain
\begin{equation}\label{eq:adp:ineq_to_rewrite}
\begin{aligned}
	&V^{adp}_{i-1}(\boldsymbol{z}) \leq l(\boldsymbol{z}, \boldsymbol{u}) + \gamma V^{adp}_{i}(\boldsymbol{A}\boldsymbol{z} + \boldsymbol{B}\boldsymbol{u}),\\
	&\forall \boldsymbol{z} \in \mathbb{R}^{6}\times \{-1, 0,1\},\;\forall \boldsymbol{u}\in \mathcal{U}(\boldsymbol{z}),\; i=1,\dots,M.
	\end{aligned}
\end{equation}
From~\cite{Wang:2014em}, we can rewrite~\eqref{eq:adp:ineq_to_rewrite} as a quadratic form
\begin{equation}\label{eq:adp:quad_form_wang}
	\begin{aligned}
	&\begin{bmatrix}\boldsymbol{z}\\ 1\end{bmatrix}^\top \boldsymbol{M}_{i}(\boldsymbol{u})\begin{bmatrix}\boldsymbol{z}\\ 1\end{bmatrix}\geq 0,\; \forall \boldsymbol{z}\in\mathbb{R}^{6}\times \{-1, 0,1\},\\
		&\;\forall \boldsymbol{u}\in \mathcal{U}(\boldsymbol{z}),\; i=1,\dots,M,
	\end{aligned}
\end{equation}
where
\begin{equation}\label{eq:adp:equation_M}
 \boldsymbol{M}_{i}(\boldsymbol{u}) = \boldsymbol{L} + \gamma \boldsymbol{G}_{i}(\boldsymbol{u}) - \boldsymbol{S}_{i-1} \in \mathbb{S}^{n_x}.
\end{equation}
is a symmetric matrix.
The matrices $\boldsymbol{S}_{i-1}, \boldsymbol{L}$ and $\boldsymbol{G}_i(\boldsymbol{u})$ are defined in Appendix~\ref{app:adp_mat}.

By noting that the state vector $\boldsymbol{z}$ includes two parts which can take only a finite set of values --- the normalized desired frequency fixed to $1$ and the previous physical input $\boldsymbol{u}_{sw}(k-1)\in \{-1,0,1\}$ --- we can explicitly enumerate part of the state-space and rewrite the quadratic form~\eqref{eq:adp:quad_form_wang} more compactly as
\begin{equation}\label{eq:adp:decomp_z}
	\begin{bmatrix}\tilde{\boldsymbol{z}}\\ 1\end{bmatrix}^\top \tilde{\boldsymbol{M}}_{i}(\boldsymbol{m})\begin{bmatrix}\tilde{\boldsymbol{z}}\\ 1\end{bmatrix}\geq 0,\; \forall \tilde{\boldsymbol{z}}\in\mathbb{R}^{8},\;\forall \boldsymbol{m}\in \mathcal{M},\; i=1,\dots,M,
\end{equation}
where $\tilde{\boldsymbol{z}}$ is the state vector without the desired frequency and $\boldsymbol{u}_{sw}(k-1)$. Moreover, $\boldsymbol{m} := \left(\boldsymbol{u}_{sw},\boldsymbol{u}_{sw,pr}\right)\in \mathcal{M}$ are all the possible combinations of current and previous switch positions satisfying the switching and integrality constraints~\eqref{eq:dmpc:constraints_u}. The detailed derivation of $\tilde{\boldsymbol{M}}(\boldsymbol{m})\in \mathbb{S}^{9}$ can be found in Appendix~\ref{app:adp_mat}.

Using the non-negativity condition of quadratic forms~\cite{Vandenberghe:1996fy}, it is easy to see that~\eqref{eq:adp:decomp_z} holds if and only if $\tilde{\boldsymbol{M}}_{i}(\boldsymbol{m})$ is positive semidefinite. Hence, 
problem~\eqref{eq:adp:underestim_first} can finally be rewritten as the following SDP
\begin{equation}\label{eq:adp:final_sdp}
 \begin{aligned}
&\text{maximize} & & \Tr\left(\boldsymbol{P}_0 \boldsymbol{\Sigma}_c\right) + 2\boldsymbol{q}_{0}^\top \boldsymbol{\mu}_c + r_0\\
&\text{subject to} & & \tilde{\boldsymbol{M}}_{i}(\boldsymbol{m}) \succeq 0, \quad \forall \boldsymbol{m}\in \mathcal{M}, \quad i=1,\dots,M\\
&&& V^{adp}_0=V^{adp}_M\\
&&& \boldsymbol{P}_i \in \mathbb{S}^{n_x},\;\boldsymbol{q}_i\in \mathbb{R}^{n_x},\; r_i \in \mathbb{R},\; i = 0,\dots,M,
 \end{aligned}
\end{equation}
which can be solved efficiently using a standard SDP solver, e.g.~\cite{mosek}. Once we obtain the solution to~\eqref{eq:adp:final_sdp}, we can define the infinite horizon tail cost to be used in problem~\eqref{eq:compmpc:problem_full} as
\begin{equation}\label{eq:adp:final_tail_cost}
	V^{adp}(\boldsymbol{z}) = \boldsymbol{z}^\top \boldsymbol{P}_0\boldsymbol{z} + 2 \boldsymbol{q}_0^\top\boldsymbol{z} + r_0.
\end{equation}

\subsection{Optimization Problem in Vector Form }\label{MPCVEC:sec}
Since we consider short horizons, we adopt a condensed MPC formulation of problem~\eqref{eq:compmpc:problem_full} with only input variables, producing a purely integer program. In this way all the possible discrete input combinations can be evaluated directly. With a sparse formulation including the continuous states within the variables, it would be necessary to solve a mixed-integer program requiring more complex computations.

Let us define the input sequence over the horizon $N$ starting at time $0$ as
\begin{equation}\label{eq:dmpc:uvec}
	\boldsymbol{U} = \begin{bmatrix} \boldsymbol{u}(0)^\top & \boldsymbol{u}(1)^\top & \dots & \boldsymbol{u}(N-1)^\top\end{bmatrix}^\top,
\end{equation}
where we have dropped the time index from $\boldsymbol{U}$ to simplify the notation. With straightforward algebraic manipulations outlined in Appendix~\ref{app:mpc_dense_mat}, it is possible to rewrite problem~\eqref{eq:compmpc:problem_full} as a parametric integer quadratic program in the initial state~$\boldsymbol{x}_0$:
\begin{equation}\label{eq:mpcvec:complete}
 \begin{aligned}
&\text{minimize} & & \boldsymbol{U}^\top\boldsymbol{Q}\boldsymbol{U} + 2\boldsymbol{f}\left(\boldsymbol{x}_0\right)^\top\boldsymbol{U}\\
&\text{subject to} & & \boldsymbol{A}_{ineq} \boldsymbol{U} \leq \boldsymbol{b}_{ineq}(\boldsymbol{x}_0)\\
&&& \mathcal{G}\boldsymbol{U} \in \{-1, 0, 1\}^{3N}.
 \end{aligned}
\end{equation}





\section{Framework for Performance Evaluation}
\label{sec:Framework for Performance Evaluation}

To benchmark our algorithm we consider a neutral point clamped voltage source inverter connected to a medium-voltage induction machine and a constant mechanical load. We consider the same model as in \cite{Geyer:ij}: a $\unit[3.3]{kV}$ and $\unit[50]{Hz}$ squirrel-cage induction machine rated at $\unit[2]{MVA}$ with a total leakage inductance of $\unit[0.25]{pu}$. On the inverter side, we assume the dc-link voltage $V_{dc}=\unit[5.2]{kV}$  to be constant and the potential of the neutral point to be fixed. The base quantities of the per unit (pu) system are the following: $V_b = \sqrt{2/3}V_{rat} = \unit[2694]{V},\; I_b = \sqrt{2}I_{rat} = \unit[503.5]{A}$ and $f_b = f_{rat} = \unit[50]{Hz}$. Quantities $V_{rat},\;I_{rat}$ and $f_{rat}$ refer to the rated voltage, current and frequency respectively. The detailed parameters are provided in Table~\ref{tab:perf:parameters}.
The switching frequency is typically in the range between $200$ and $\unit[350]{Hz}$ for medium-voltage inverters~\cite{Geyer:ij}.
If not otherwise stated, all simulations were done at rated torque, nominal speed and fundamental frequency of $\unit[50]{Hz}$.

\begin{table}[h]
\caption{Rated values and parameters of the drive \cite{Geyer:ij}}
\centering
\begin{adjustbox}{max width=\columnwidth}
  \begin{tabular}{@{}lc lc @{}c lc@{}}
\toprule
\multicolumn{5}{c}{Induction Motor} & \multicolumn{2}{c}{Inverter}\\
\cmidrule{1-4}  \cmidrule{6-7}
Voltage & $\unit[3300]{V}$ & $R_s$ & \unit[0.0108]{pu} & & $V_{dc}$ & $\unit[1.930]{pu}$\\
Current & $\unit[356]{A}$ & $R_r$ & \unit[0.0091]{pu} & & $x_{c}$ & $\unit[11.769]{pu}$\\
Real power & $\unit[1.587]{MW}$ & $X_{ls}$ & \unit[0.1493]{pu} & &  & \\
Apparent power & $\unit[2.035]{MVA}$ & $X_{lr}$ & \unit[0.1104]{pu} & &  & \\
Frequency & $\unit[50]{Hz}$ & $X_{m}$ & \unit[2.3489]{pu} & &  & \\
Rotational speed & $\unit[596]{rpm}$ & &  & &  & \\
     \bottomrule
  \end{tabular}
\end{adjustbox}
  \label{tab:perf:parameters}
\end{table}

We consider an idealized model with the semiconductors switching instantaneously.  As such, we neglect second-order effects like deadtimes, controller delays, measurement noise, observer errors, saturation of the machine's magnetic material, variations of the parameters and so on. 
This is motivated by the fact that, using a similar model, previous simulations~\cite{Geyer:2009fu} showed a very close match with the experimental results in~\cite{Papafotiou:2009gs}.
All the steady-state simulations in the following sections were also performed with model mismatch of $\pm \unit[1]{\%}$ in all the parameters of Table~\ref{tab:perf:parameters} showing negligible variations in the THD. However, we omit these benchmarks since an exhaustive sensitivity analysis is out of the scope of this paper.



\section{Achievable Performance in Steady State}
\label{sec:Achievable Performance in Steady State}

We performed closed loop simulations in steady-state operation in MATLAB to benchmark the achievable performance in terms of THD and switching frequency. The system was simulated for $4$ periods before recording to ensure it reaches steady-state operation. The THD and switching frequency were computed over simulations of $20$ periods. The discount factor was chosen as $\gamma = 0.95$ and the switching frequency filter parameters as $r_1=r_2=800$ in order to get a smooth estimate. The weighting $\delta$ was chosen such that the switching frequency is around $\unit[300]{Hz}$. The infinite horizon estimation SDP \eqref{eq:adp:final_sdp} is formulated using YALMIP \cite{YALMIP} with $M=50$ Bellman iterations and solved offline using MOSEK \cite{mosek}. Note that in case of a change in the systems parameters, e.g.\ the dc-link voltage or the rotor speed, the tail cost has to be recomputed. However, it possible to precompute offline and store several quadratic tail costs for different possible parameters and evaluate the desired one online without significant increase complexity.

For comparison, we simulated the drive system also with the direct MPC controller described in~\cite{Geyer:ij} (denoted DMPC) tuned in order to have the same switching frequency by adjusting the weighting parameter $\lambda_u$.

The integer optimization problems were solved using Gurobi Optimizer~\cite{gurobi}. Numerical results with both approaches are presented in Table~\ref{tab:achievable_performance:simulation_results}. Note that the choice of the solver does not influence the THD or the switching frequency and we would have obtained the same results with another optimization software.

%
%
%

\begin{table}[h]
\caption{Simulation Results with ADP and with DMPC from~\cite{Geyer:fva} at switching frequency $\unit[300]{Hz}$}
\centering
  \begin{tabular}{@{}l cc@{} c cc@{}}
\toprule
&\multicolumn{2}{c}{ADP} & & \multicolumn{2}{c}{DMPC~\cite{Geyer:fva}}\\
\cmidrule{2-3}  \cmidrule{5-6}
 & $\delta$  & $\unit[\mathrm{THD}]{[\%]}$ & & $\lambda_u$ & $\unit[\mathrm{THD}]{[\%]}$\\
 \cmidrule{2-3}  \cmidrule{5-6}
$N = 1$ & $4$ & $5.24$ &  & $0.00235$ & $5.44$\\
$N = 2$ & $5.1$ & $5.13$ &  & $0.00690$ & $5.43$\\
$N = 3$ & $5.5$ & $5.10$ &  & $0.01350$ & $5.39$\\
$N = 10$ & $10$ & $4.80$ & & $0.10200$ & $5.29$\\
     \bottomrule
  \end{tabular}
  \label{tab:achievable_performance:simulation_results}
\end{table}

Our method, with a horizon of $N=1$ provided both a THD improvement over the DMPC formulation in \cite{Geyer:ij} with ${N=10}$ and a drastically better numerical speed. This showed how choosing a meaningful cost function can provide good control performance without recourse to long horizons.
Moreover, we also performed a comparison with longer horizons ${N=2}$, ${N=3}$ and ${N=10}$. Our method, with horizon ${N=10}$ would give an even greater reduction in THD to $\unit[4.80]{\%}$.


\section{FPGA Implementation}
\label{sec:FPGA Implementation}


\subsection{Hardware Setup}
\label{sub:Hardware Setup}
\begin{figure}[hbt]
\begin{center}
  \ifTwoColumn
  \includegraphics[width=0.8\columnwidth]{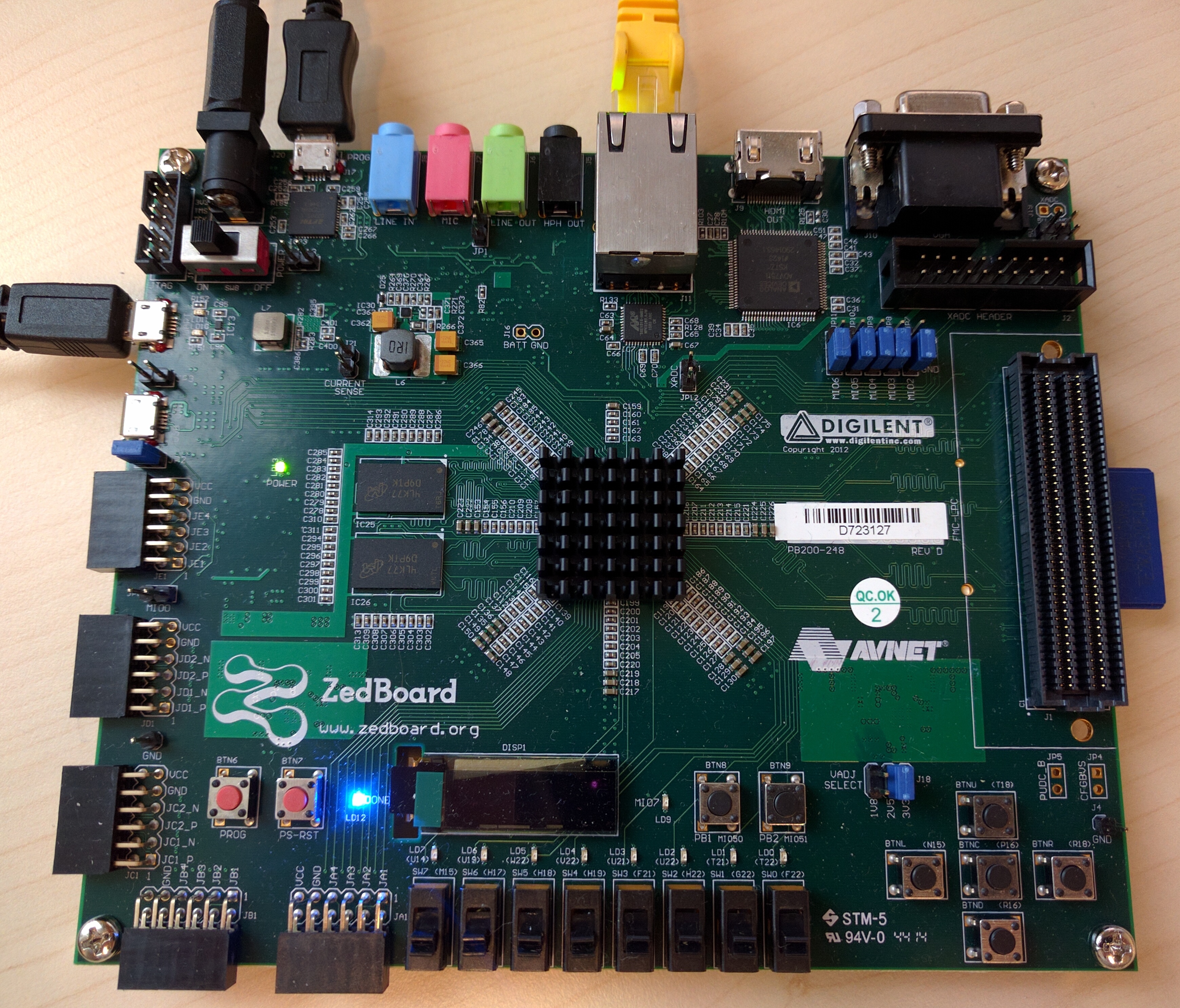}
  \else
  \includegraphics[width=0.6\textwidth]{img/zedboard}
  \fi
  \caption{Zedboard Evaluation Board used for HIL Tests. The controller runs on the FPGA while the plant is simulate on the laptop. The states and input vectors are passed via the ethernet cable (yellow). The micro-usb cable on the left side provides a UART interface with the laptop used to print if there are any problems in the communication. The cable in the top left corner is connected to the power supply while the other micro-usb cable next to it provides access to the {USB-JTAG} interface to program the FPGA module.}
  \label{img:fpga_implementation:zedboard}
  \end{center}
\end{figure}

We implemented the control algorithm on a Xilinx Zynq (xc7z020) \cite{Anonymous:7K9ZD1ex}, a low cost FPGA, running at approximately $\unit[150]{MHz}$ mounted on the Zedboard evaluation module\cite{Anonymous:L9mDSuzn}; see Figure~\ref{img:fpga_implementation:zedboard}.
The control algorithm was coded in C++ using the PROTOIP Toolbox \cite{Suardi:2015ia}. The FPGA vendor’s HLS tool Xilinx Vivado HLS~\cite{XilinxInc:ZYWTdbK2} was used to convert the written code to VHDL defining the Programmable Logic connections.

\subsection{Algorithm Description}
\label{sub:Algorithm Description}

We now present a detailed description of how the controller within the dotted lines in Figure~\ref{fig:block_diagram} was implemented on the FPGA.

The updates in OSC and FLT were implemented as simple matrix multiplications. The solver for the integer problem~\eqref{eq:mpcvec:complete} was implemented with a simple exhaustive search algorithm for three reasons: first, the tail cost approximation provides good performance with very few horizon steps while considering a relatively small number of input combinations; second, the structure of the problem allows us to evaluate both the inequalities and cost function for multiple input sequences in parallel; third, the FPGA logic is particularly suited for highly pipelined and/or parallelized operations, which are at the core of exhaustive search.

To exploit the FPGA architecture, we implemented our algorithm in fixed-point arithmetic using custom data types defined in Vivado HLS~\cite{XilinxInc:ZYWTdbK2}. In particular, we used  $4$ integer and $0$ fractional bits to describe the integer inputs  and $2$ integer and $22$ fractional bits to describe the states and the cost function values. This choice is given by the minimum number of bits necessary to describe these quantities from floating-point simulations in Section~\ref{sec:Achievable Performance in Steady State}. Note that the exhaustive search algorithm does not suffer from any accumulation of rounding error because it consists entirely of independent function evaluations, in contrast to iterative optimization algorithms~\cite{nocedal2006numerical}.

\begin{algorithm}

  \ifTwoColumn
  \else
          \ifTechReport
          \else
          \setstretch{1} 
          \fi
  \fi
\caption{Controller Algorithm}
\label{alg:fpga_implementation:controller_algorithm}
\begin{algorithmic}[1]
\Function{ComputeMPCinput}{$T^*(k), \boldsymbol{x}_{ph}(k)$}

  \Statex \textit{Data:} ${\boldsymbol{x}_{osc}(k-1)}$, ${\boldsymbol{x}_{sw}(k-1)}$, ${\boldsymbol{p}(k-1)}$ and ${\boldsymbol{u}_{sw}(k-1)}$

  \Statex \textit{Parameters:} $\boldsymbol{U}^{seq} \in \mathbb{Z}^{6\times 27^N}$,  $J_{ub}\in \mathbb{R}$

  \Statex \textit{Initialize:} $\boldsymbol{J} \in \mathbb{R}^{27^N}$, $J_{min}\in \mathbb{R}$ and $i_{min}\in \mathbb{N}$

  \Statex
  \Statex \textit{Execute Filter and Oscillator to Obtain Initial State:}
  \If{change in $T(k)$} \label{alg:fpga_implementation:controller_algorithm:filters_beg}
  \State $\boldsymbol{x}_{osc}(k) \gets$ Reset according to~\eqref{eq:drive_system:torque}
  \Else
  \State $\boldsymbol{x}_{osc}(k) \gets \boldsymbol{A}_{osc}\boldsymbol{x}_{osc}(k-1)$
  \EndIf

  \State $\boldsymbol{x}_{sw}(k) \gets \boldsymbol{A}_{sw}\boldsymbol{x}_{sw}(k-1) + \boldsymbol{B}_{sw}\boldsymbol{p}(k-1)$

  \State $\boldsymbol{x}_0 \gets \begin{bmatrix}\boldsymbol{x}_{ph}(k)^\top & \boldsymbol{x}_{osc}(k)^\top &  \boldsymbol{x}_{sw}(k)^\top & \boldsymbol{u}_{sw}(k-1)^\top \end{bmatrix}^\top$\label{alg:fpga_implementation:controller_algorithm:filters_end}
  \Statex

  \Statex \textit{Precompute Vectors:}

  \State $\boldsymbol{f}(\boldsymbol{x}_0) \gets$ Compute from~\eqref{eq:mpcvec:fx0} \label{alg:fpga_implementation:controller_algorithm:prec_fx0}
  \State $\boldsymbol{b}_{ineq}(\boldsymbol{x}_0) \gets$ Compute from~\eqref{eq:mpcvec:ineq1},~\eqref{eq:mpcvec:ineq2} \label{alg:fpga_implementation:controller_algorithm:prec_b_ineq}
  \Statex

  \Statex \textit{Loop 1 - Compute Cost Function Values:}
  \For{$i=1,\dots,27^N$}\label{alg:fpga_implementation:controller_algorithm:cost_fun_beg}

  \State $\boldsymbol{u} \gets \boldsymbol{U}^{seq}_{(:,i)}$
  \State $\boldsymbol{u}_{(4:6)} \gets \boldsymbol{p}(k) = \| \boldsymbol{u}_{(1:3)} - \boldsymbol{u}_{sw}(k-1)\|_1$ \label{alg:fpga_implementation:controller_algorithm:line_p}

  \If{$\boldsymbol{A}_{ineq}\boldsymbol{u}\leq \boldsymbol{b}_{ineq}(\boldsymbol{x}_0)$}

  \State $\boldsymbol{J}_{(i)} \gets  \boldsymbol{u}^\top \boldsymbol{Q}\boldsymbol{u}+ 2\boldsymbol{f}(\boldsymbol{x}_0)^\top \boldsymbol{u}$ \label{alg:fpga_implementation:controller_algorithm:store_J1}

  \Else

  \State $\boldsymbol{J}_{(i)} \gets J_{ub}$ \label{alg:fpga_implementation:controller_algorithm:store_J2}

  \EndIf
  \EndFor \label{alg:fpga_implementation:controller_algorithm:cost_fun_end}

  \Statex
  \Statex \textit{Loop 2 - Find Minimum:}
  \State $J_{min} \gets J_{ub}$, $i_{min} \gets 1$ \label{alg:fpga_implementation:controller_algorithm:minimum_beg}

  \For{$i=1,\dots,27^N$}
    \If{$\boldsymbol{J}_{(i)} \leq J_{min}$}
      \State $J_{min} \gets \boldsymbol{J}_{(i)}$
      \State $i_{min} \gets i$
    \EndIf
  \EndFor \label{alg:fpga_implementation:controller_algorithm:minimum_end}

 \Statex
 \Statex \textit{Return Results}
 \State $\boldsymbol{u}_{sw}(k)\gets \boldsymbol{U}^{seq}_{(1:3,i_{min})}$ and $\boldsymbol{p}(k)\gets \boldsymbol{U}^{seq}_{(4:6,i_{min})}$
 \State \Return  $\boldsymbol{u}_{sw}(k)$

\EndFunction

\end{algorithmic}
\end{algorithm}

We provide pseudo-code for our method in Algorithm~\ref{alg:fpga_implementation:controller_algorithm}. From Figure~\ref{fig:block_diagram}, the controller receives the required torque $T^*(k)$ and the motor states $\boldsymbol{x}_{ph}(k)$ and returns the switch positions $\boldsymbol{u}_{sw}(k)$. From line~\ref{alg:fpga_implementation:controller_algorithm:filters_beg} to line~\ref{alg:fpga_implementation:controller_algorithm:filters_end} the oscillator OSC and the filter FLT are updated to compute the new initial state $\boldsymbol{x}_0$ for the optimization algorithm. Note that if there is a change in the required torque then the oscillator states $\boldsymbol{x}_{osc}(k)$ are reset to match the new $T^*(k)$.
Line~\ref{alg:fpga_implementation:controller_algorithm:prec_fx0} and~\ref{alg:fpga_implementation:controller_algorithm:prec_b_ineq} precompute the vectors in problem~\eqref{eq:mpcvec:complete} depending on $\boldsymbol{x}_0$.

The main loop iterating over  all input combinations is split into two subloops: Loop 1, which is completely decoupled and can be parallelized; and Loop 2, which can only be pipelined.

Loop 1 from line~\ref{alg:fpga_implementation:controller_algorithm:cost_fun_beg} to~\ref{alg:fpga_implementation:controller_algorithm:cost_fun_end} computes the cost function values for every combination $i$ and stores it into vector $\boldsymbol{J}$. All the possible input sequence combinations are saved in the static matrix $\boldsymbol{U}^{seq}$. For every loop cycle,  sequence $i$ is saved into variable $\boldsymbol{u}$. Then, in line~\ref{alg:fpga_implementation:controller_algorithm:line_p}, the value of $\boldsymbol{p}(k)$ is updated inside $\boldsymbol{u}$ with $\boldsymbol{u}_{sw}(k-1)$ according to~\eqref{eq:dmpc:rel_p_u1} and~\eqref{eq:dmpc:rel_p_u2}.
If $\boldsymbol{u}$ satisfies the constraint $\boldsymbol{A}_{ineq}\boldsymbol{u}\leq\boldsymbol{b}_{ineq}(\boldsymbol{x}_0)$, then the cost function is stored in $\boldsymbol{J}_{(i)}$ (line \ref{alg:fpga_implementation:controller_algorithm:store_J1}). Otherwise $\boldsymbol{J}_{(i)}$ is set to a high value $J_{ub}$. Note that, even though it would bring considerable speed improvements, we do not precompute offline the quadratic part $\boldsymbol{u}^\top \boldsymbol{Q}\boldsymbol{u}$ of the cost and the left side of the inequality $\boldsymbol{A}_{ineq}\boldsymbol{u}$ since it would also require enumeration over inputs at the previous control cycle used in line~\ref{alg:fpga_implementation:controller_algorithm:line_p}.
 %
 %
 %

 Each iteration of this loop is independent from the others and can therefore be parallelized efficiently.

Loop 2 from line~\ref{alg:fpga_implementation:controller_algorithm:minimum_beg} to~\ref{alg:fpga_implementation:controller_algorithm:minimum_end} is a simple loop iterating over the computed cost function values to find the minimum and save it into $J_{min}$. Every iteration depends sequentially on $J_{min}$ which is accessed and can be modified at every $i$. Thus, in this form it is not possible to parallelize this loop, although it can be pipelined.

\subsection{Circuit Generation}
\label{sub:Circuit Generation}

In Vivado HLS~\cite{XilinxInc:ZYWTdbK2} it is possible to specify directives to optimize the circuit synthesis according to the resources available on the target board. Loop 1 and Loop 2 were pipelined and the preprocessing operations from line~\ref{alg:fpga_implementation:controller_algorithm:filters_beg} to~\ref{alg:fpga_implementation:controller_algorithm:prec_b_ineq} parallelized. We generated the circuit for the algorithm~\ref{alg:fpga_implementation:controller_algorithm} with horizons $N=1$ and $N=2$ at frequency $\unit[150]{MHz}$ (clock cycle of $\unit[7]{ns}$). The resources usage and the timing estimates are displayed in Table~\ref{tab:fpga_implementation:resources}. Since timing constraints were met, there was no need to parallelize Loop 1 to reduce computation time.

\begin{table}[h]
\caption{Resources Usage and Timing Estimates for Implementation on the Xilinx Zynq FPGA (xc7z020) running at $\unit[150]{MHz}$}
\centering
\begin{tabular}{@{}l l l l@{}}
\toprule
& & $N=1$ & $N=2$\\
\cmidrule{1-4}
\multirow{4}{*}{FPGA Resources} & LUT & $15127\; (\unit[28]{\%})$ & $31028\; (\unit[58]{\%})$\\
 & FF & $11156\; (\unit[10]{\%})$ & $20263\; (\unit[19]{\%})$\\
 & BRAM & $6\; (\unit[2]{\%})$& $21\; (\unit[7]{\%})$\\
 & DSP & $89\; (\unit[40]{\%})$& $201\; (\unit[91]{\%})$\\
\cmidrule{1-4}
 Clock Cycles & & 371 & 1953\\
 Delay & & $\unit[2.60]{\mu s}$& $\unit[13.67]{\mu s}$\\
\bottomrule
\end{tabular}
  \label{tab:fpga_implementation:resources}
\end{table}

Note that for $N=2$ we are using already $\unit[91]{\%}$ of the DSP multipliers. This is due to the limited amount of resources available on our chosen low-cost hardware.


\section{Hardware In The Loop Tests}
\label{sec:Hardware In The Loop Tests}
We performed hardware in the loop (HIL) experiments using the controller FPGA fixed-point implementation developed in Section~\ref{sec:FPGA Implementation} and the machine model described in Section~\ref{sec:Framework for Performance Evaluation}.

The control loop was operated using the PROTOIP toolbox~\cite{Suardi:2015ia}: the plant model was simulated on a Macbook Pro 2.8 GHz Intel Core i7 with 16GB of RAM while the control algorithm was entirely executed on the Zedboard development board described in Section~\ref{sub:Hardware Setup}.

\subsection{Steady State}
\label{sub:Steady State}

The controller was benchmarked in HIL in steady-state operation to compare its performance to the achievable performance results obtained in Table~\ref{tab:achievable_performance:simulation_results}. We chose the same controller parameters as in Section~\ref{sec:Achievable Performance in Steady State}.

%
%

\begin{figure*}\centering
	\subfloat[Three-phase stator currents (solid lines) with their references (dashed lines).]{\label{fig:hil_tests:steady_results:currents}\includegraphics[width=0.325\textwidth]{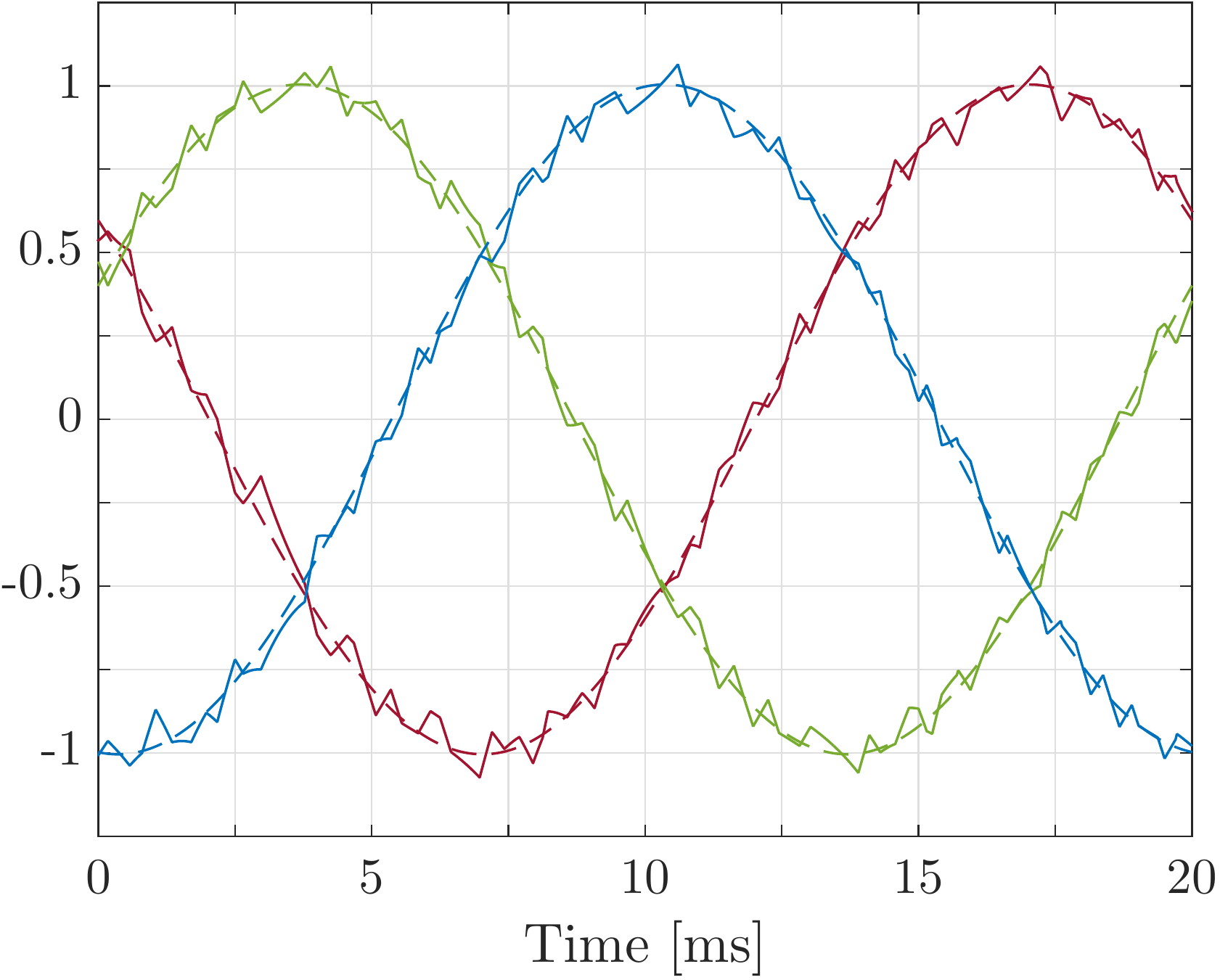}}\hfill
	\subfloat[Stator current spectrum.]{\label{fig:hil_tests:steady_results:spectrum}\includegraphics[width=0.337\textwidth]{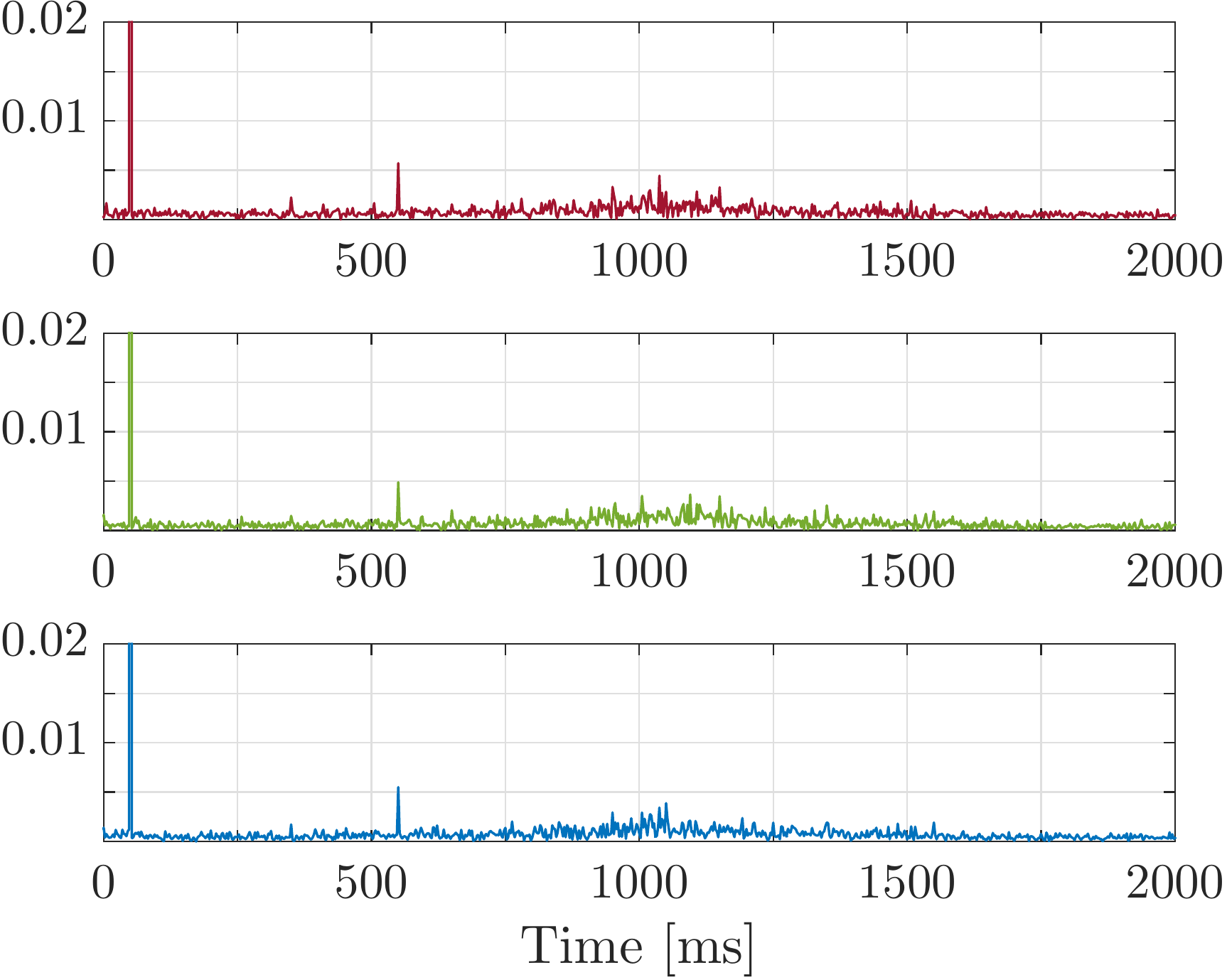}}\hfill
	\subfloat[Three-phase switch positions inputs.]{\label{fig:hil_tests:steady_results:inputs}\includegraphics[width=0.318\textwidth]{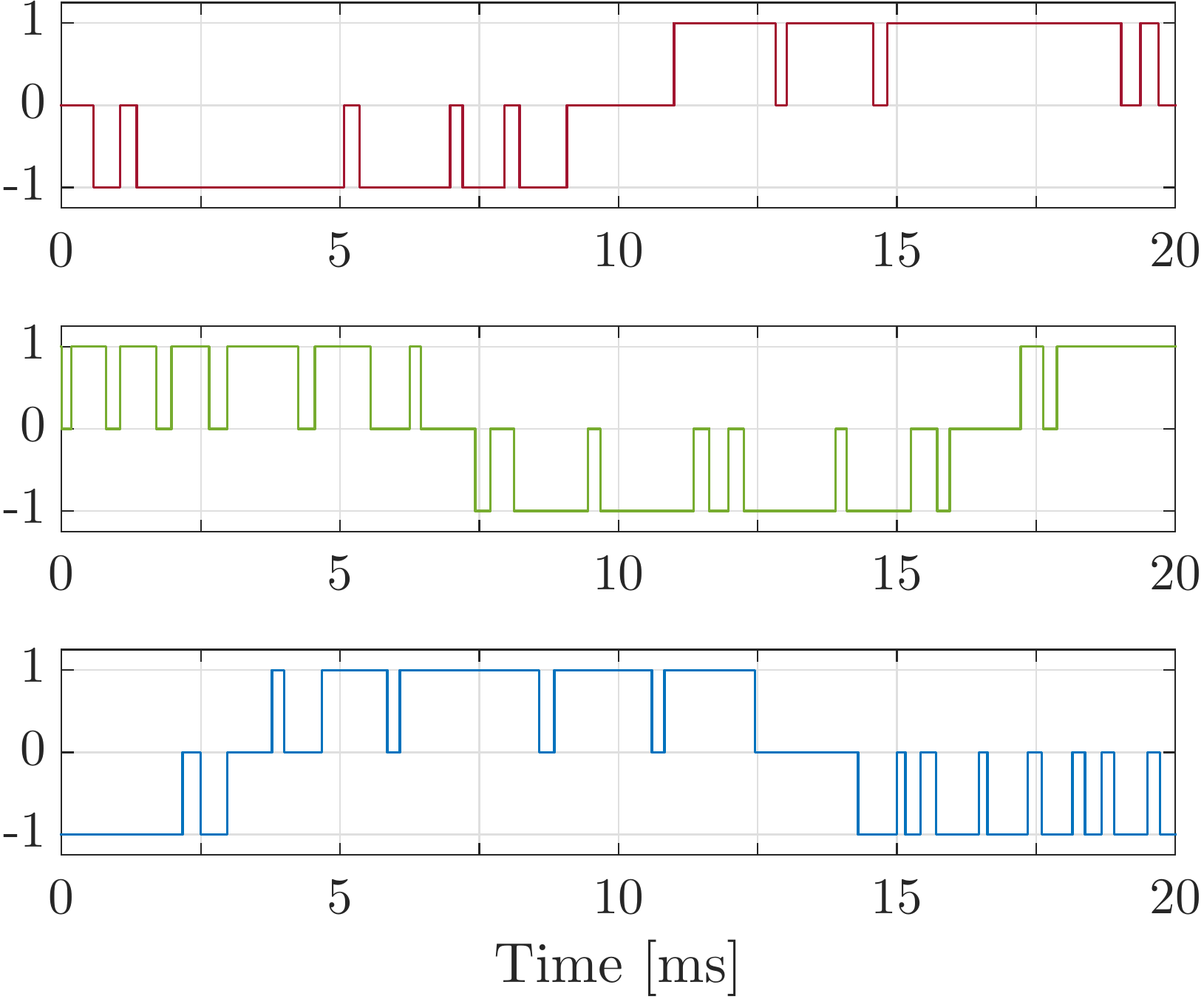}}
	\caption{Waveforms produced during HIL Tests by the direct model predictive controller at steady state operation, at full speed and rated torque. Horizon of $N=1$ is used. The switching frequency is approximately $\unit[300]{Hz}$ and the current THD is $\unit[5.23]{\%}$.}
	\label{fig:hil_tests:steady_results}
\end{figure*}

The HIL tests for horizon $N=1$ are shown in Figure~\ref{fig:hil_tests:steady_results} in the per unit system. The three-phase stator currents are displayed over a fundamental period in Figure~\ref{fig:hil_tests:steady_results:currents}, the three spectra are shown in Figure~\ref{fig:hil_tests:steady_results:spectrum} with THD of $\unit[5.23]{\%}$ and the input sequences are plotted in Figure~\ref{fig:hil_tests:steady_results:inputs}.

From the experimental benchmarks with horizon $N=1$ and $N=2$ we obtained ${\mathrm{THD} =\unit[5.23]{\%}}$ and ${\mathrm{THD}=\unit[5.14]{\%}}$ respectively. As expected, these results are very close to the simulated ones in Table~\ref{tab:achievable_performance:simulation_results}. The slight difference ($\sim \unit[0.01]{\%}$) comes from the fixed-point implementation of the oscillator OSC and the filter FLT in Figure~\ref{fig:block_diagram}.

\subsection{Transients}
\label{sub:Transients}

\begin{figure*}\centering
\subfloat[Three-phase stator currents (solid lines) with their references (dashed lines).]{\label{fig:hil_tests:trans_results:currents}\includegraphics[width=0.32\textwidth]{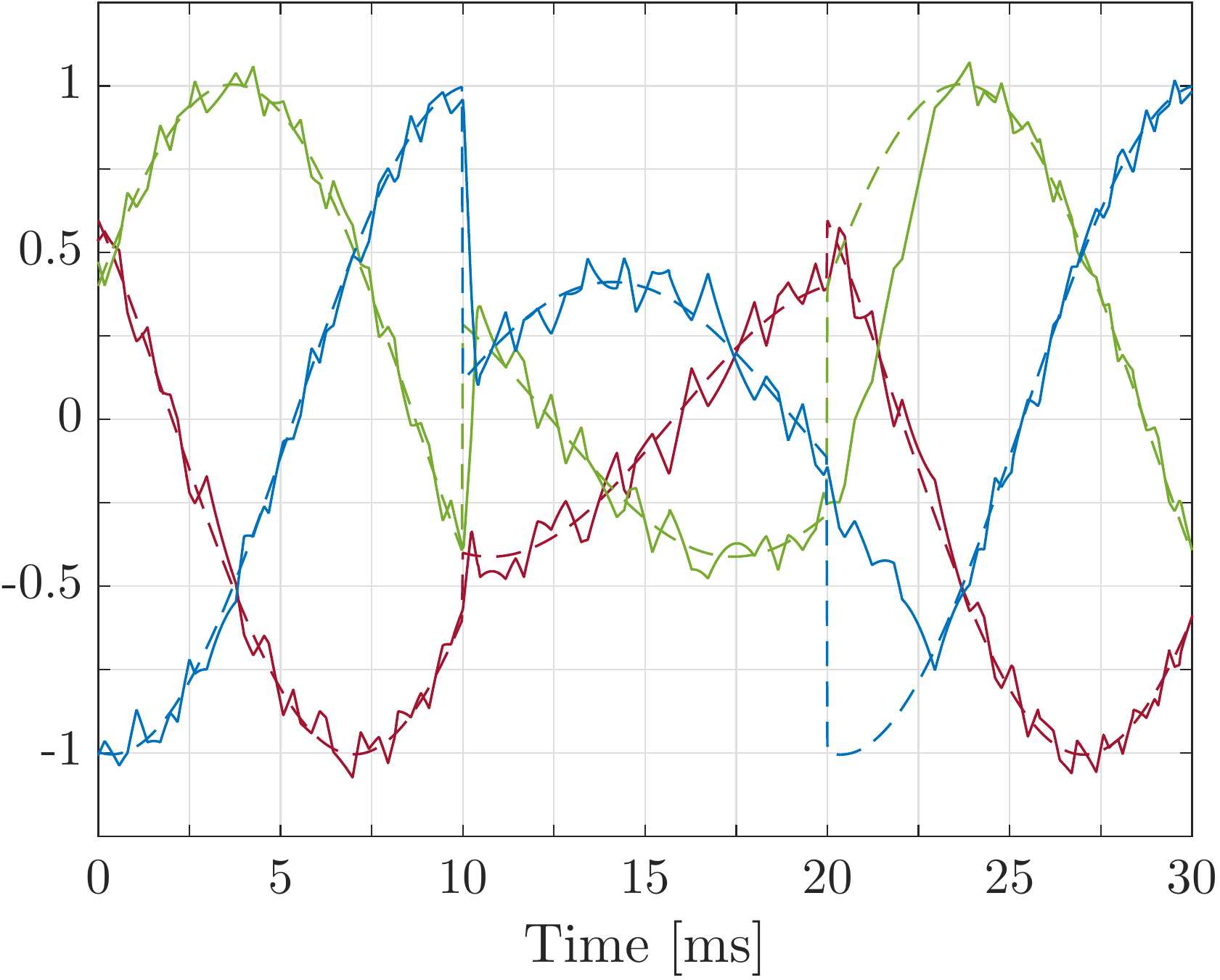}}\hfill
\subfloat[Actual (solid line) and reference (dashed line) torques.]{\label{fig:hil_tests:trans_results:torque}\includegraphics[width=0.32\textwidth]{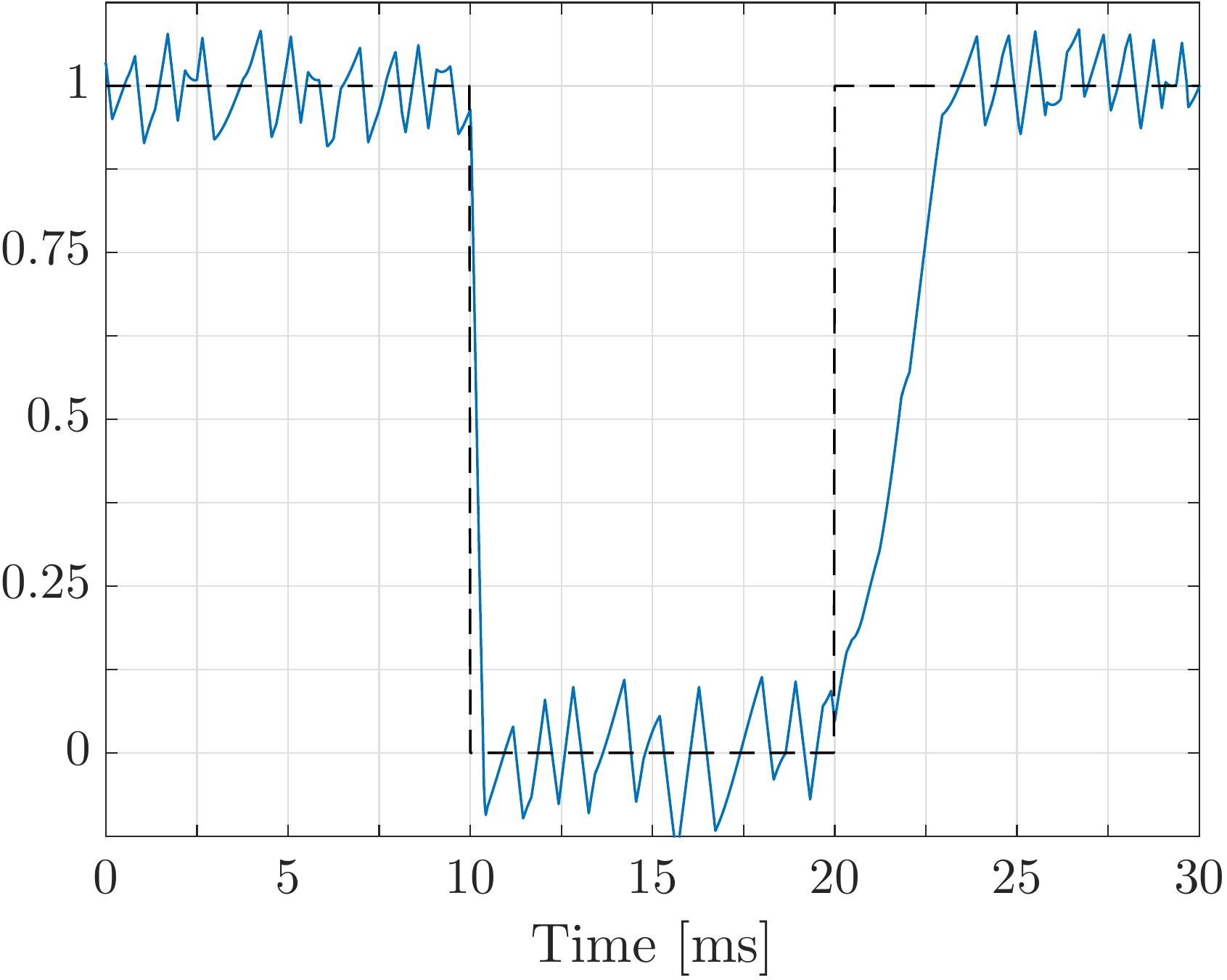}}\hfill
\subfloat[Three-phase switch positions inputs.]{\label{fig:hil_tests:trans_results:inputs}\includegraphics[width=0.305\textwidth]{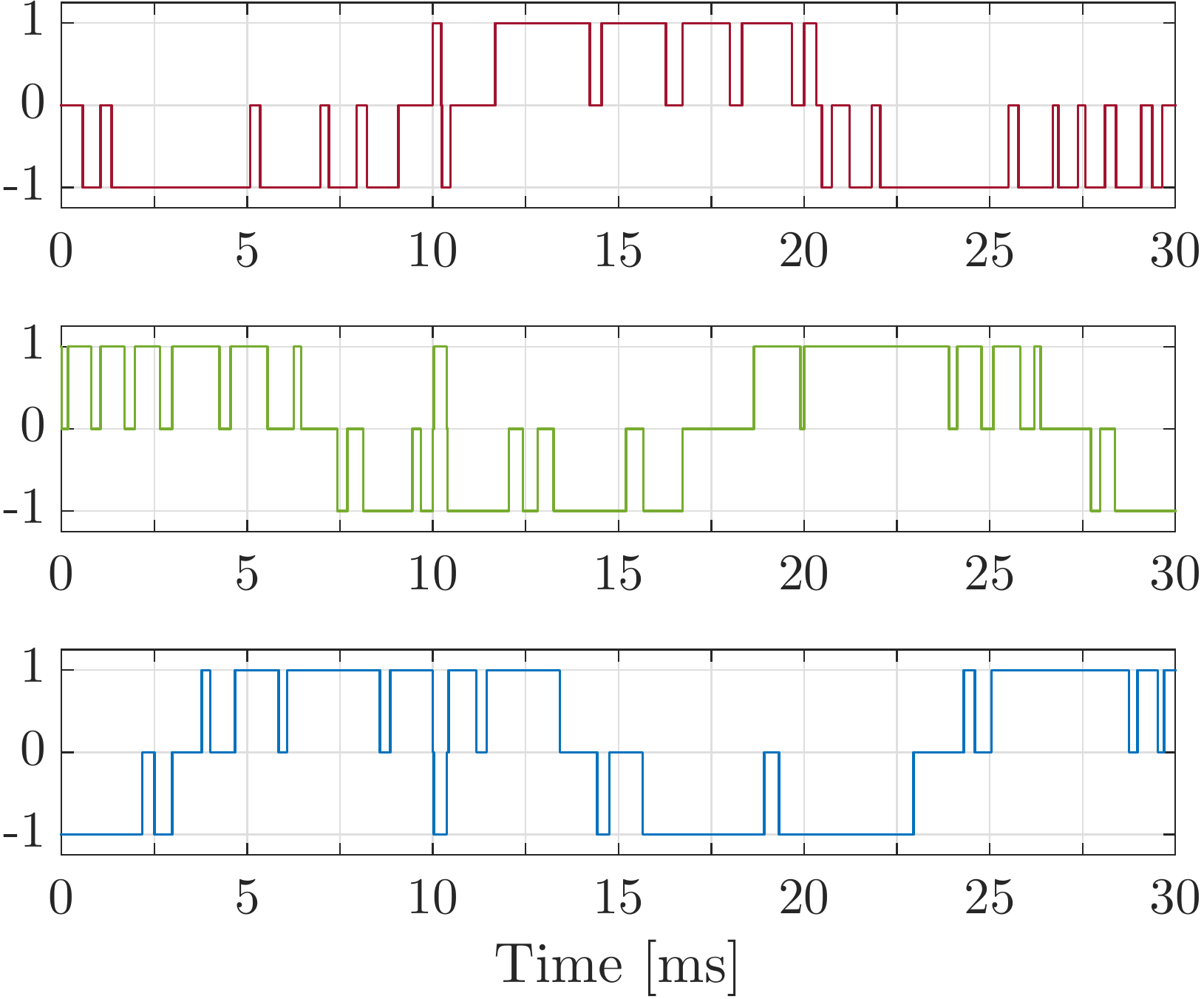}}
  \caption{Reference torque steps produced by the direct model predictive controller in HIL tests with horizon $N=1$.}
  \label{fig:hil_tests:trans_results}
\end{figure*}

One of the main advantages of direct MPC is the
fast transient response \cite{Geyer:ij}. We tested torque transients in HIL with the same tuning parameters as in the steady state benchmarks. At nominal speed,
reference torque steps were imposed; see
Figure~\ref{fig:hil_tests:trans_results:torque}. These steps were
translated into different current references to track, as shown in
Figure~\ref{fig:hil_tests:trans_results:currents}, while the computed inputs
are shown in Figure~\ref{fig:hil_tests:trans_results:inputs}.

The torque step from $1$ to $0$ in the per unit system presented an extremely short settling
time of $\unit[0.35]{ms}$ similar to deadbeat control approaches~\cite{rodriguez2012predictive}. This was achieved by inverting the voltage
applied to the load. Since we prohibited switchings between $-1$ and
$1$ in~\eqref{eq:dmpc:constraints_u:fsw1}
and~\eqref{eq:dmpc:constraints_u:fsw2}, the voltage inversion was
performed in $2T_s$ via an intermediate zero switching position.

Switching from $0$ to $1$ torque produced much slower response time of approximately $\unit[3.5]{ms}$. This was due to the limited
available voltage in the three-phase admissible switching positions. As shown in Figure~\ref{fig:hil_tests:trans_results:inputs},
during the second step at time $\unit[20]{ms}$, the phases $b$ and $c$ saturated at the values $+1$ and $-1$ respectively for the
majority of the transient providing the maximum available voltage that
could steer the currents to the desired values.



These results match the simulations of the DMPC formulation in~\cite{Geyer:ij} in terms of settling time showing that our method possesses the fast dynamical behavior during transients typical of direct current MPC.

As noted in~\cite{Geyer:ij}, having a longer horizon or a better predictive behavior does not significantly improve the settling times. This is because the benefit of longer prediction obtainable by extending the horizon or adopting a powerful final stage costs is reduced by the saturation of the inputs during the transients.


\begin{figure}[h]
	\begin{center}
    \begin{tikzpicture}
      \centering
      \begin{axis}[
        width=\columnwidth,
        x tick label style={/pgf/number format/.cd,%
        set thousands separator={}, 
        },
        xbar,
        enlarge y limits=0.5,
        xmin = 0, xmax = 1050,
        y = 4.26em,
        bar width=1.136em,
        axis y line*=left,
        legend style={at={(1,0.90)},
        },
        xlabel = {Execution Time $[\mu s]$},
        symbolic y coords={$N=2$,$N=1$},
        ytick=data,
        nodes near coords,
        nodes near coords align=horizontal,
        axis on top,
         extra x ticks={25},  
         extra x tick label = {$T_s$},
         extra x tick style={    
        xticklabel pos=right,   
         tick style={draw = none},
         major grid style={black!50, densely dashed},
         xmajorgrids=true            
         }
        ]
        \addplot[every node near coord/.append style={anchor = east, color = redCol}, area legend,
                draw = redCol, fill = redCol!10,
                error bars/error bar style={black},
                error bars/.cd,
                x dir=plus, x explicit] coordinates {(621.20,$N=1$)+- (119.98,0.0) (750.40,$N=2$)+- (216.15,0.0)};
        \addplot[every node near coord/.append style={xshift = 0.3em, color = blueCol}, draw = blueCol, fill = blueCol!30, area legend] coordinates {(5.76,$N=1$) (17.27,$N=2$)};
        \legend{Gurobi, FPGA}
      \end{axis}
    \end{tikzpicture}
	\end{center}
	\caption{Execution times required by the Xilinx Zynq FPGA (xc7z020) to execute our controller based on an ADP formulation (blue) and using Gurobi Optimizer~\cite{gurobi} to solve the formulation in~\cite{Geyer:ij} on a Macbook Pro with Intel Core i7 2.8 GHz and 16GB of RAM}
	\label{fig:hil_tests:timing}
\end{figure}
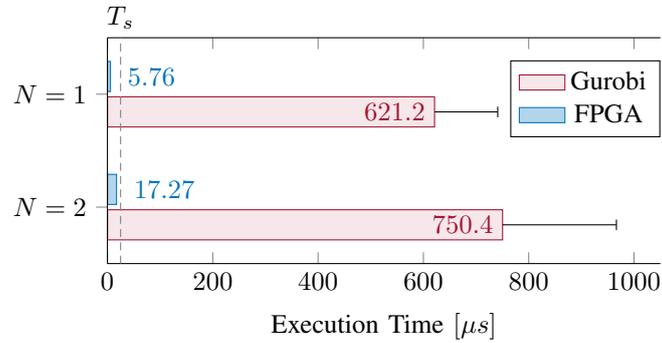

\subsection{Execution Time}
\label{sub:Execution Time}

To show that the controller is able to run on cheap hardware within $T_s = \unit[25]{\mu s}$, we measured the time the FPGA took to execute Algorithm~\ref{alg:fpga_implementation:controller_algorithm} for horizon $N=1$ and $N=2$. Since there are no available DMPC sphere decoding algorithm execution times, we compared our results to the time needed to solve the DMPC formulation in~\cite{Geyer:fva} for the same horizon lengths on a Macbook Pro with Intel Core i7 2.8 GHz and 16GB of RAM using the commercial integer program solver Gurobi Optimizer~\cite{gurobi} which implements an efficient Branch-and-Bound algorithm. The results are shown in Figure~\ref{fig:hil_tests:timing}.

The FPGA execution times were $\unit[5.76]{\mu s}$ and $\unit[17.27]{\mu s}$ for horizon $N=1$ and $N=2$ respectively. Note that they presented a slight overhead of approximately $\unit[3.5]{\mu s}$ compared to the estimates in Table~\ref{tab:fpga_implementation:resources} since the measured times included the time needed to exchange the input-output data from the FPGA to the ARM processor through the RAM memory. Without the overhead, the estimated FPGA computing times obtained by the circuit generation are exact; see~\cite{XilinxInc:ZYWTdbK2}.

Note that the time needed by the FPGA to compute the control algorithm is deterministic with zero variance. This makes our HIL implementation particularly suited for real-time applications. Furthermore, it is important to underline that the method we propose is the \emph{only} method available that can produce integer optimal solutions to this problem achieving this performance in $\unit[25]{\mu s}$ sampling time.

The execution times needed by Gurobi optimizer were ${\unit[621.2 \pm 119.98]{\mu s}}$ and $\unit[750.40 \pm 216.15]{\mu s}$ for horizons $N=1$ and $N=2$ respectively. The non-negligible standard deviation appeared because of the branch-and-bound algorithm implemented in Gurobi. However, since we are considering real-time applications, we are interested in the worst case number of visited nodes which is always the whole tree of combinations, i.e. $27^N$. Note that the DMPC formulation was solved in~\cite{Geyer:ij} using a different branch-and-bound algorithm based on the sphere decoding algorithm~\cite{Hassibi:2005gu}, but the worst case number of visited nodes cannot be easily reduced because of the $\mathcal{NP}-$hardness of the problem.



\section{Conclusion and Future Work}\label{sec:conclusions}
This work proposes a new computationally efficient direct model
predictive control (MPC) scheme for current reference tracking in
power converters. We extended the problem formulation in~\cite{Geyer:ij}
and~\cite{Geyer:fva} in order to include a switching frequency
estimator in the system state and rewrite the optimal control problem
as a regulation one. To reduce the horizon length and decrease the
computational burden while preserving good control performances, we
estimated the infinite horizon tail cost of the MPC problem formulation using approximate dynamic programming (ADP).

Steady-state simulation results show that with our method requiring
short horizons, it is possible to obtain better performance than the direct MPC formulation in \cite{Geyer:ij} with long
horizons. This is due to the predictive
behavior of the tail cost function obtained with ADP.

The control algorithm was implemented in fixed-point arithmetic on the low size Xilinx Zynq FPGA (xc7z020) for horizons $N=1$ and $N=2$. Hardware in the loop (HIL) tests during steady-state operation showed an almost identical performance to the simulation results. We also
performed transient simulations where our
proposed approach exhibited the same very fast dynamic response
as the direct MPC described in~\cite{Geyer:ij}. Moreover, we showed that our algorithm can run within the sampling time of $\unit[25]{\mu s}$ by measuring the execution time on the FPGA. Results showed that only $\unit[5.76]{\mu s}$ and $\unit[17.27]{\mu s}$ are required to run our controller for horizons $N=1$ and $N=2$ respectively.


Direct MPC can also be applied to more complex schemes such as modular multilevel converters (MMC)~\cite{Riar:2014kn}. While it is possible to derive a complete MMC model that could be used in an MPC approach, the number of switching levels per horizon stage exponentially increases with the number of converter levels. As stated in~\cite{Geyer:fva}, long-horizon predictive power is expected to be even more beneficial with MMCs. We believe that our method, making use of short prediction horizons and long predictions using an approximate value function could be  applied effectively to MMCs with more levels because it is still possible to evaluate on commercially available FPGAs the multilevel feasible switching combinations over very short horizons within the required sampling time.

There are several future directions to be investigated. Given the system  design there are several symmetries in the model  that could be exploited to increase the controller horizon without requiring more computational power. Regarding the frequency estimation, other filters with different orders could be implemented and their parameters chosen optimally by solving an optimization problem instead of performing manual tuning. Moreover, it would be interesting to benchmark other ADP tail cost functions (e.g. higher order polynomials) to understand which ones best approximate the infinite horizon tail cost and produce the best overall control performance.
\section*{Acknowledgment}
The authors would like to thank Andrea Suardi and Bulat Khusainov for their advice regarding FPGA implementation, in particular the PROTOIP toolbox.


%

\appendices%

\section{Physical System Matrices}\label{app:phys_sys_mat}
The matrices corresponding to the continuous-time physical system in \eqref{eq:complete_phys:continuous_time} are
\begin{equation*}
\boldsymbol{D} = \begin{bmatrix}
 -\dfrac{\strut 1}{\strut \tau_s} & 0 & \dfrac{\strut X_m}{\strut \tau_r D} & \omega_r\dfrac{\strut X_m}{\strut D}\\
 0 & -\dfrac{\strut 1}{\strut \tau_s} & -\omega_r \dfrac{\strut X_m}{\strut D} & \dfrac{\strut X_m}{\strut \tau_r D}\\
 \dfrac{\strut X_m}{\strut \tau_r} & 0 & -\dfrac{\strut 1}{\strut \tau_r} & -\omega_r\\
 0 & \dfrac{\strut X_m}{\strut \tau_r} & \omega_r & -\dfrac{\strut 1}{\strut \tau_r}
 \end{bmatrix}
\end{equation*}
\begin{equation*}
\boldsymbol{E} = \dfrac{X_r}{D}\dfrac{V_{dc}}{2}\begin{bmatrix}
1 & 0\\
0 & 1\\
0 & 0\\
0 & 0
\end{bmatrix}\boldsymbol{P}, \qquad \boldsymbol{F} = \begin{bmatrix}
 1 & 0 & 0 & 0\\
 0 & 1 & 0 & 0
 \end{bmatrix}.
\end{equation*}

\section{ADP Formulation}\label{app:adp_mat}
The matrices defining the quadratic form are
\begin{equation*}
	\begin{aligned}
		\boldsymbol{S}_{i-1} \coloneqq &\begin{bmatrix}\boldsymbol{P}_{i-1} & \boldsymbol{q}_{i-1}\\ \boldsymbol{q}_{i-1}^{\top} & r_{i-1}\end{bmatrix}, \quad \boldsymbol{L} \coloneqq \begin{bmatrix}\boldsymbol{C}^\top\boldsymbol{C} & \boldsymbol{0}_{n_x \times 1}\\ \boldsymbol{0}_{n_x \times 1}^\top & 0\end{bmatrix}\\
		&\boldsymbol{G}_{i}(\boldsymbol{u}) \coloneqq \begin{bmatrix} \boldsymbol{\Psi}^{(i)} & \boldsymbol{\Phi}^{(i)}(\boldsymbol{u}) \\ \boldsymbol{\Phi}^{(i)}(\boldsymbol{u})^\top & \Gamma ^{(i)}(\boldsymbol{u})\end{bmatrix},
	\end{aligned}
\end{equation*}
with
\begin{equation*}
	\begin{aligned}
		\boldsymbol{\Psi}^{(i)} &= \boldsymbol{A}^\top\boldsymbol{P}_{i-1}\boldsymbol{A}\\
		\boldsymbol{\Phi}^{(i)}(\boldsymbol{u}) &= \boldsymbol{A}^\top\boldsymbol{P}_{i}\boldsymbol{B}\boldsymbol{u} + \boldsymbol{A}^\top\boldsymbol{q}_{i}\\
		\Gamma ^{(i)}(\boldsymbol{u}) &= \boldsymbol{u}^\top\boldsymbol{B}^\top\boldsymbol{P}_i\boldsymbol{B}\boldsymbol{u} + 2\boldsymbol{q}_i^\top\boldsymbol{B}\boldsymbol{u} + r_i,
	\end{aligned}
\end{equation*}
for $i=1,\dots,M$.

The quadratic form decomposition can be derived as follows. For every ${\boldsymbol{m} = \left(\boldsymbol{u}_{sw},\boldsymbol{u}_{sw, pr}\right)\in \mathcal{M}}$, we can define the input
\begin{equation*}
	\boldsymbol{u}_{\boldsymbol{m}} = \begin{bmatrix}\boldsymbol{u}_{sw}\\
	\left\|\boldsymbol{u}_{sw}- \boldsymbol{u}_{sw,pr}\right\|_1
 \end{bmatrix},
\end{equation*}
and the matrix $\boldsymbol{M}_{i}(\boldsymbol{u}_{\boldsymbol{m}})$ using \eqref{eq:adp:equation_M}. Now we can decompose the vector in quadratic form \eqref{eq:adp:quad_form_wang} using the state definition~\eqref{eq:dmpc:state_def}
\begin{equation*}
	\begin{bmatrix}
	\boldsymbol{z}^\top & 1
	\end{bmatrix}^\top =\begin{bmatrix}
	\boldsymbol{z}_{ph}^\top & \boldsymbol{z}_{osc}^{\top}  & 	\boldsymbol{z}_{sw, 1:2}^{\top} & \boldsymbol{z}_{sw, 3}^{\top} & \boldsymbol{z}_{upr}^{\top}  & 1
	\end{bmatrix}^\top.
\end{equation*}
The matrix $\boldsymbol{M}_{i}(\boldsymbol{u}_{\boldsymbol{m}})$  can also be decomposed in the same fashion into smaller block matrices as follows
\begin{equation*}
	\begin{bmatrix}\boldsymbol{z}_{ph}\\ \boldsymbol{z}_{osc}\\ \boldsymbol{z}_{sw, 1:2}\\ \boldsymbol{z}_{sw,3} \\\boldsymbol{z}_{upr} \\  1\end{bmatrix}^\top
		\left[\begin{array}{ccc|ccc}
		 \multicolumn{3}{c|}{\multirow{3}{*}{$\boldsymbol{M}_{i,11}$}} & \multirow{3}{*}{$\boldsymbol{M}_{i,12}$} & \multirow{3}{*}{$\boldsymbol{M}_{i,13}$}&  \multirow{3}{*}{$\boldsymbol{M}_{i,14}$}\\
		  &  &   & \\
		  &  &   & \\
		 \hline
		 \multicolumn{3}{c|}{\boldsymbol{M}_{i,12}^\top} & \boldsymbol{M}_{i,22} &\boldsymbol{M}_{i,23}  & \boldsymbol{M}_{i,24}\\
 		 \multicolumn{3}{c|}{\boldsymbol{M}_{i,13}^\top} & \boldsymbol{M}_{i,23}^\top &\boldsymbol{M}_{i,33}  & \boldsymbol{M}_{i,34}\\
		 \multicolumn{3}{c|}{\boldsymbol{M}_{i,14}^\top} & \boldsymbol{M}_{i,24}^\top &\boldsymbol{M}_{i,34}^\top  & \boldsymbol{M}_{i,44}
		\end{array}\right]
	\begin{bmatrix}\boldsymbol{z}_{ph}\\ \boldsymbol{z}_{osc}\\ \boldsymbol{z}_{sw, 1:2}\\ \boldsymbol{z}_{sw,3} \\\boldsymbol{z}_{upr} \\  1\end{bmatrix}\geq 0,
\end{equation*}
where the dependency $\boldsymbol{M}_{i, **}(\boldsymbol{u}_{m}),\; \boldsymbol{m}\in \mathcal{M}$ has been neglected to simplify the notation. The first row and first column block matrices have the first and second dimensions respectively equal to the length of vector $\begin{bmatrix} \boldsymbol{z}_{ph}^\top & \boldsymbol{z}_{osc}^\top & \boldsymbol{z}_{sw,1:2}^\top\end{bmatrix}^\top$. Since $\boldsymbol{z}_{upr}= \boldsymbol{u}_{sw, pr}$ and $\boldsymbol{z}_{sw,3} = 1$ (normalized desired switching frequency), we can rewrite the quadratic form as
\begin{equation}\label{eq:app:adp_quad_form_condensed}
	\begin{bmatrix}\boldsymbol{z}_{ph}\\ \boldsymbol{z}_{osc} \\ \boldsymbol{z}_{sw, 1:2}\\ 1\end{bmatrix}^\top
		\left[\begin{array}{ccc|c}
		 \multicolumn{3}{c|}{\multirow{3}{*}{$\boldsymbol{M}_{i,11}$}} & \multirow{3}{*}{$\boldsymbol{\Psi}_{i,1}$} \\
		  &  & & \\
      &  & & \\
		 \hline
		 \multicolumn{3}{c|}{\boldsymbol{\Psi}_{i,1}^\top} &  \boldsymbol{\Psi}_{i,2}\\
		\end{array}\right]
	\begin{bmatrix}\boldsymbol{z}_{ph}\\ \boldsymbol{z}_{osc} \\ \boldsymbol{z}_{sw, 1:2}\\ 1\end{bmatrix}\geq 0,
\end{equation}
where
\begin{equation*}
\begin{aligned}
	\boldsymbol{\Psi}_{i,1} &= \boldsymbol{M}_{i,13}\boldsymbol{z}_{upr} + \boldsymbol{M}_{i,12} + \boldsymbol{M}_{i,14} \\
	\boldsymbol{\Psi}_{i,2} &= \boldsymbol{z}_{upr}^\top\boldsymbol{M}_{i,33}\boldsymbol{z}_{upr} +   2\boldsymbol{M}_{i,23}\boldsymbol{z}_{upr}  + 2\boldsymbol{M}_{i,34}^\top \boldsymbol{z}_{upr} \\
	&\quad + 2\boldsymbol{M}_{i,24} + \boldsymbol{M}_{i, 22}
\end{aligned}
\end{equation*}
Therefore, we will denote the matrix in \eqref{eq:app:adp_quad_form_condensed} as $\tilde{\boldsymbol{M}}_{i}(\boldsymbol{m})$ and the quadratic form vectors as $\begin{bmatrix}\tilde{\boldsymbol{z}}^\top & 1\end{bmatrix}^\top$.

\section{Dense Formulation of the MPC Problem}\label{app:mpc_dense_mat}
By considering the input sequence \eqref{eq:dmpc:uvec} and the state sequence over the horizon  denoted as
\begin{equation*}
	\boldsymbol{X} = \begin{bmatrix} \boldsymbol{x}(0)^\top & \boldsymbol{x}(1)^\top & \dots & \boldsymbol{x}(N)^\top\end{bmatrix}^\top,
\end{equation*}
the system dynamics \eqref{eq:compmpc:problem_full:dyn} with initial state constraint \eqref{eq:compmpc:problem_full:initst}  can be written as
\begin{equation}\label{eq:mpcvec:state_dynamics}
	\boldsymbol{X} = \mathcal{A}\boldsymbol{x}_0 + \mathcal{B}\boldsymbol{U},
\end{equation}
where $\mathcal{A}$ and $\mathcal{B}$ are
\begin{equation*}
\begin{aligned}
	\mathcal{A} \coloneqq \begin{bmatrix}
 \boldsymbol{I}\\
 \boldsymbol{A}\\
 \vdots\\
 \boldsymbol{A}^{N}
 \end{bmatrix}
\end{aligned}, \quad \mathcal{B} \coloneqq \begin{bmatrix}
 	\boldsymbol{0} &  \multicolumn{2}{c}{\dots} &\boldsymbol{0}\\
 	\boldsymbol{B} &  &  &  \multirow{2}{*}{\vdots} \\
 	\boldsymbol{A}\boldsymbol{B} & \boldsymbol{B} & & \\
 	\vdots &  \ddots & \ddots & \boldsymbol{0} \\
 	\boldsymbol{A}^{N-1}\boldsymbol{B} & \dotsm & \boldsymbol{A}\boldsymbol{B} &\boldsymbol{B}
 \end{bmatrix}.
\end{equation*}
Let us separate cost function \eqref{eq:compmpc:problem_full:cost_fun} in two parts: the cost from stage $0$ to $N-1$ and the tail cost. The former can be rewritten as
\begin{equation}\label{eq:mpcvec:cost_fun1}
	\begin{aligned}
		&\sum_{k=0}^{N-1}\gamma^{k}\|\boldsymbol{C}\boldsymbol{x}(k)\|_2^2 = \boldsymbol{X}^\top \mathcal{H} \boldsymbol{X} \\
		& = \boldsymbol{U}^\top \mathcal{B}^\top \mathcal{H}\mathcal{B}\boldsymbol{U} + 2\left(\mathcal{B}^\top\mathcal{H}\mathcal{A}\boldsymbol{x}_0\right)^\top \boldsymbol{U} + \mathrm{const}(\boldsymbol{x}_0),
	\end{aligned}
\end{equation}
where the last equality is obtained by plugging in \eqref{eq:mpcvec:state_dynamics} and the term $\mathrm{const}(\boldsymbol{x}_0)$ is a constant depending on the initial state. Matrix $\mathcal{H}$ is defined as
\begin{equation}
	\mathcal{H} = \begin{bmatrix}
 	\boldsymbol{C}^\top \boldsymbol{C} & & & & \boldsymbol{0}\\
 	 & \gamma \boldsymbol{C}^\top \boldsymbol{C}& & & \multirow{3}{*}{\vdots}\\
 	 & & \ddots & & \\
 	 & &  & \gamma^{N-1} \boldsymbol{C}^\top \boldsymbol{C} & \\
 	 \boldsymbol{0} & \multicolumn{3}{c}{\dotsm}  & \boldsymbol{0}
 \end{bmatrix}.
\end{equation}
In order derive the tail cost, let us write the last stage as
\begin{equation}\label{MPCVEC:final_state}
	\boldsymbol{x}(N) = \boldsymbol{A}^{N}\boldsymbol{x}_0 + \mathcal{B}_{end}\boldsymbol{U},
\end{equation}
where $\mathcal{B}_{end}$ is the last row of $\mathcal{B}$ used to compute the last state. Using equations \eqref{MPCVEC:final_state} and \eqref{eq:adp:final_tail_cost}, the tail cost can be rewritten as
\begin{equation}\label{MPCVEC:cost_fun2}
	\begin{aligned}
		&V^{adp}(\boldsymbol{x}(N)) = \boldsymbol{x}(N)^\top \boldsymbol{P}_0 \boldsymbol{x}(N) + 2\boldsymbol{q}_0^\top \boldsymbol{x}(N) + r_0\\
		&= \boldsymbol{U}^\top \left(\mathcal{B}_{end}^\top \boldsymbol{P}_0\boldsymbol{B}_{end}\right)\boldsymbol{U} + 2\left(\mathcal{B}_{end}^\top\boldsymbol{P}_0\boldsymbol{A}^{N}\boldsymbol{x}_0 + \mathcal{B}_{end}^\top\boldsymbol{q}_0\right)^\top\boldsymbol{U}\\
		&\quad +\mathrm{const}(\boldsymbol{x}_0).
	\end{aligned}
\end{equation}

By combining \eqref{eq:mpcvec:cost_fun1} and \eqref{MPCVEC:cost_fun2} according to \eqref{eq:compmpc:problem_full:cost_fun}, we obtain the full cost function
\begin{equation*}
	\boldsymbol{J} = \boldsymbol{U}^\top\boldsymbol{Q}\boldsymbol{U} + 2\boldsymbol{f}\left(\boldsymbol{x}_0\right)^\top\boldsymbol{U} + \mathrm{const}(\boldsymbol{x}_0),
\end{equation*}
with
\begin{align}
		\boldsymbol{Q} &= \mathcal{B}^\top\mathcal{H}\mathcal{B} + \gamma^{N}\mathcal{B}_{end}^\top \boldsymbol{P}_0\mathcal{B}_{end} \label{eq:mpcvec:Q}\\
		\boldsymbol{f}\left(\boldsymbol{x}_0\right) & = \left(\mathcal{B}^\top\mathcal{H}\mathcal{A} + \gamma^{N}\mathcal{B}_{end}\boldsymbol{P}_0\mathcal{A}^{N}\right)\boldsymbol{x}_0 \label{eq:mpcvec:fx0}\\
		&\quad + \gamma^{N}\mathcal{B}_{end}^\top \boldsymbol{q}_0 \nonumber.
	\end{align}

We now rewrite the constraints of problem~\eqref{eq:compmpc:problem_full} in vector form. Inequalities~\eqref{eq:dmpc:constraints_u:fsw1} with $k=0,\dots,N-1$ can be written as
\begin{equation}\label{eq:mpcvec:ineq1}
\mathcal{R} \boldsymbol{U} \leq \mathcal{S}\boldsymbol{X} \iff (\mathcal{R} - \mathcal{S}\mathcal{B})\boldsymbol{U}\leq \mathcal{S}\mathcal{A}\boldsymbol{x}_0,
\end{equation}
where in the term on the right we substituted \eqref{eq:mpcvec:state_dynamics}.

Similarly, constraint \eqref{eq:dmpc:constraints_u:fsw2} with $k=0,\dots,N-1$ can be written as
\begin{equation}\label{eq:mpcvec:ineq2}
	\begin{aligned}
		\left\|\boldsymbol{T}\boldsymbol{u}(k)\right\|_\infty \leq 1 \iff \mathcal{F}\boldsymbol{U} \leq \boldsymbol{1},
	\end{aligned}
\end{equation}
where $\boldsymbol{1}$ is a vector of ones of appropriate dimensions. Matrices $\mathcal{R}, \mathcal{S}$ and $\mathcal{F}$ are
\begin{equation*}
	\mathcal{R} = \begin{bmatrix}
	\boldsymbol{I} - \boldsymbol{T} & &\\
	& \ddots &\\
	& & \boldsymbol{I} - \boldsymbol{T}\\
	-\boldsymbol{I} - \boldsymbol{T} & &\\
	& \ddots &\\
	& & -\boldsymbol{I} - \boldsymbol{T}
 \end{bmatrix},
\end{equation*}
\begin{equation*}
 \mathcal{S} = \begin{bmatrix}
 	\boldsymbol{W} & & & \boldsymbol{0}\\
 	 & \ddots & & \multirow{4}{*}{\vdots}\\
 	 &  & \boldsymbol{W} & \\
  	-\boldsymbol{W} & & & \\
 	 & \ddots & & \\
 	 &  & -\boldsymbol{W} & \boldsymbol{0}
 \end{bmatrix},
\;\mathcal{F} = \begin{bmatrix}
	\boldsymbol{T} & &\\
     & \ddots &\\
	 & & 	\boldsymbol{T}\\
	-\boldsymbol{T} & &\\
     & \ddots &\\
	 & & 	-\boldsymbol{T}
\end{bmatrix}.
\end{equation*}

Let us define matrix $\mathcal{G}$ extrapolating all the switch positions $\boldsymbol{u}_{sw}(k)$ from $\boldsymbol{U}$ as
\begin{equation*}
	 \mathcal{G} = \begin{bmatrix}
 	 G & & \\
 	  & \ddots &\\
 	  & & G
  \end{bmatrix}.
\end{equation*}

We can now merge \eqref{eq:mpcvec:ineq1} and \eqref{eq:mpcvec:ineq2} into a
single inequality $\boldsymbol{A}_{ineq} \boldsymbol{U} \leq
\boldsymbol{b}_{ineq}(\boldsymbol{x}_0)$ and rewrite
\eqref{eq:compmpc:problem_full} neglecting the constant terms in the cost
function obtaining the result in \eqref{eq:mpcvec:complete}.


\balance



\ifCLASSOPTIONcaptionsoff
  \newpage
\fi



\bibliographystyle{IEEEtran}
\bibliography{bibliography}
\end{document}